\documentclass[11pt]{article}

\usepackage[T1]{fontenc}
\usepackage{lmodern}
\usepackage{amsmath,amssymb,amsthm,mathtools}
\usepackage{mathrsfs}
\usepackage{xcolor}
\usepackage{enumitem}
\usepackage{young}
\usepackage{tikz}
\usetikzlibrary{arrows.meta}
\usepackage[margin=1in]{geometry}
\usepackage[colorlinks=true,linkcolor=blue!55!black,urlcolor=blue!55!black]{hyperref}

\newtheorem{theorem}{Theorem}[section]
\newtheorem{proposition}[theorem]{Proposition}
\newtheorem{lemma}[theorem]{Lemma}
\newtheorem{corollary}[theorem]{Corollary}
\newtheorem{remark}[theorem]{Remark}
\theoremstyle{definition}
\newtheorem{definition}[theorem]{Definition}
\newtheorem{example}[theorem]{Example}

\newcommand{\A}{\mathcal A}
\newcommand{\Hh}{\mathcal H}
\newcommand{\I}{I}
\newcommand{\FPF}{\mathrm{FPF}}
\newcommand{\PrB}{\widehat P_{\mathrm{rB}}}
\newcommand{\PcB}{\widehat P_{\mathrm{cB}}}
\newcommand{\shape}{\operatorname{shape}}
\newcommand{\Des}{\operatorname{Des}}
\newcommand{\Asc}{\operatorname{Asc}}
\newcommand{\Fix}{\operatorname{Fix}}
\newcommand{\SYT}{\operatorname{SYT}}

\newcommand{\wt}{\operatorname{wt}}
\newcommand{\topJ}{\operatorname{top}_{J}}
\newcommand{\tr}{\mathsf t}
\newcommand{\Oinf}{\mathscr O_{\infty}}
\newcommand{\mcell}{\mathbf m}
\newcommand{\ncell}{\mathbf n}
\newcommand{\mGamma}{\Gamma^{\mathsf{row}}}
\newcommand{\nGamma}{\Gamma^{\mathsf{col}}}
\newcommand{\mAsc}{\operatorname{Asc}^{\mathsf{row}}}
\newcommand{\nAsc}{\operatorname{Asc}^{\mathsf{col}}}
\newcommand{\RSinsert}{\xleftarrow{\,\mathsf{RS}\,}}
\newcommand{\rBinsert}{\xleftarrow{\,\mathsf{rB}\,}}
\newcommand{\cBinsert}{\xleftarrow{\,\mathsf{cB}\,}}
\newcommand{\rowtableau}[1]{\vcenter{\hbox{\begin{Young}#1\end{Young}}}}

\title{Levi restriction and cells in the Gelfand \texorpdfstring{$S_n$}{S(n)}-graphs}
\author{Yifeng Zhang
\\
School of Mathematical Sciences \\
South China Normal University \\
{\tt calvinz314159@gmail.com}}
\date{}

\begin{document}
\maketitle

\begin{abstract}
In the row and column Gelfand $W$-graphs for the symmetric group, cells
are the strongly connected components, and molecules are the connected
components obtained by retaining only bidirected edges.
Marberg and Zhang classified the molecules by modified Beissinger
insertion shape and conjectured that every molecule is a cell.
We prove this conjecture for both reduced graphs by showing that insertion
shapes weakly decrease in dominance order along coefficient arrows and
hence are constant on each cell.
We obtain this monotonicity by identifying each insertion shape with the
transpose of the corresponding Levi top and using the Levi filtration
to show that these tops weakly increase along arrows.
\end{abstract}

\section{Introduction}

Marberg and Zhang constructed two Gelfand $W$-graphs for the symmetric
group $S_d$, the row graph $\mGamma$ and the column graph $\nGamma$
\cite{MarbergZhangGelfand2022}.  Each gives a model for the Hecke algebra
of $S_d$ in which every irreducible representation over $\mathbb Q(v)$
occurs once.  As in the $W$-graphs of Kazhdan and Lusztig
\cite{KazhdanLusztig1979}, the edges record the action of the Hecke
generators on a distinguished basis.  We use the reduced graphs
throughout.  For distinct vertices $y,z$, there is an arrow
$y\longrightarrow z$ if the basis vector indexed by $z$ has nonzero
coefficient in the action of some generator $H_i$ on the vector indexed
by $y$.

A \emph{cell} is a strongly connected component of this directed graph.
Keeping only the edges for which the reverse arrow also occurs gives
connected components called \emph{molecules}.  Marberg and Zhang
classified the molecules by modified Beissinger insertion
\cite{MarbergZhangInsertion2024}.  In the row graph, two vertices lie
in the same molecule precisely when their row insertion tableaux have
the same shape.  The analogous statement holds for column insertion
in the column graph.  Proposition~\ref{prop:reduced-molecules} gives
the classification in the reduced convention used here.
Marberg and Zhang conjectured that each molecule is a cell
\cite{MarbergZhangInsertion2024}.  Since the vertices of a molecule
are mutually reachable, the question is whether a directed path can
leave one insertion shape and later return to it.

To state the result, let $\I_d=\{w\in S_d:w^2=1\}$.  We index the row
vertices by the ascending completions of $w\in\I_d$ and the column
vertices by their descending completions, as in
Subsection~\ref{subsec:reduced-graphs}.  Let $\PrB(w)$ and $\PcB(w)$
be the restricted modified row and column Beissinger tableaux, and put
\[
 \lambda_{\mathrm r}(w)=\shape\PrB(w),\qquad
 \lambda_{\mathrm c}(w)=\shape\PcB(w).
\]
For partitions of $d$, write $\alpha\preceq\beta$ if
$\sum_{i=1}^j\alpha_i\leq\sum_{i=1}^j\beta_i$ for every $j\geq1$,
taking missing parts to be zero.  A downset for this order contains
every partition dominated by one of its members.
For example, the column graph for $S_3$ has an arrow
$(12)\longrightarrow_{\ncell}(13)$ from shape $(3)$ to shape $(2,1)$,
with $(2,1)\preceq(3)$.

\begin{samepage}
\begin{theorem}\label{thm:main-intro}
Let $d\geq1$.
\begin{enumerate}[label=\textup{(\alph*)}]
\item The fibers of $\lambda_{\mathrm r}$ are the cells of $\mGamma$,
      and each is a single molecule.  Every row arrow
      $y\longrightarrow_{\mcell}z$ satisfies
      $\lambda_{\mathrm r}(z)\preceq\lambda_{\mathrm r}(y)$.
\item The fibers of $\lambda_{\mathrm c}$ are the cells of $\nGamma$,
      and each is a single molecule.  Every column arrow
      $y\longrightarrow_{\ncell}z$ satisfies
      $\lambda_{\mathrm c}(z)\preceq\lambda_{\mathrm c}(y)$.
\end{enumerate}
In each model, the $\mathbb Z[v,v^{-1}]$-span of the canonical basis
vectors whose insertion shapes belong to a given dominance downset is a
saturated Hecke submodule, so the quotient is torsion-free over
$\mathbb Z[v,v^{-1}]$.
\end{theorem}
\end{samepage}

Let $J$ index the simple roots of a Levi subalgebra.  In a highest-weight
crystal, keep only the arrows with colors in $J$.  The component
containing a vertex $b$ has a unique highest element $b^+$, characterized
by $\widetilde e_j b^+=0$ for all $j\in J$.  We call the $J$-weight of
$b^+$ the \emph{Levi top} of $b$ and denote it by $\topJ(b)$.
For the type-$A_{d-1}$ Levi modules of polynomial degree $d$ considered
below, we record this weight as a partition of $d$.

In Kashiwara's crystal theory, the multiplicity of a Levi constituent
of highest weight $\mu$ is the number of $J$-components with highest
weight $\mu$ \cite{Kashiwara1991}.
For the restriction of an irreducible module to a standard Levi,
Berenstein and Zelevinsky express these multiplicities as counts of
integer points in polyhedra \cite[\S2.4]{BerensteinZelevinsky2001}.

Lusztig showed that the highest-weight filtration of a based module
has a basis on each isotypic quotient induced by the original basis
\cite[\S27.1]{Lusztig1993}.  For complex reductive groups, Baumann,
Gaussent and Littelmann proved that the Mirkovi\'c--Vilonen basis of a
tensor product of irreducibles is $L$-perfect.  Under restriction to a
standard Levi, basis subsets span the filtration by highest weights
bounded above.  On each isotypic quotient the induced basis respects
a decomposition into irreducibles
\cite[\S3.2, Theorem~3.4]{BaumannGaussentLittelmann2024}.
In their study of separable Weyl group elements, Gossow and Yacobi use
the highest weights of restricted crystal components to define a preorder
on crystal elements \cite[Definition~4.4]{GossowYacobi2025}.

For $d\geq2$, the type-$C$ row model and the type-$D$ column model
realize each fixed-point block on the weight space $(1,\ldots,1)$ of a
Levi module.  There $H_i=E_iF_i-v^{-1}\operatorname{id}$ preserves the
span of the canonical vectors whose Levi tops dominate a given weight.
Thus $y\longrightarrow z$ implies $\topJ(y)\preceq\topJ(z)$.
A rank-two calculation connects each top fiber by bidirected edges,
and monotonicity keeps every cell in one fiber.  Hence each fiber is
both a molecule and a cell.  Comparing its two Schur characteristics
identifies its top with the transpose of its insertion shape.
Transposition reverses dominance, so shape downsets correspond to top
upsets, whose canonical vectors span saturated Hecke submodules.

This paper is organized as follows. Section~\ref{sec:preliminaries} fixes the graph and insertion conventions
and computes the Schur characteristics of the shape fibers.
Section~\ref{sec:levi-top-criterion} gives a criterion for identifying
cells from highest-weight labels and proves the restriction and
connectivity results used to apply it.  Section~\ref{sec:row-model}
constructs the row module from a type-$C$ parabolic quotient and proves
part~\textup{(a)}.  Section~\ref{sec:column-model} uses a type-$D$
parabolic quotient to prove part~\textup{(b)} for the column module.
Each of these last two sections also proves the corresponding integral
filtration statement and ends with an explicit cell for $S_6$.

%

\section*{AI Usage Announcement}
OpenAI Codex was used in revising this manuscript to draft and edit
exposition, develop parts of the proofs, and write and run programs
to check examples in small ranks.

\section{Preliminaries}
\label{sec:preliminaries}

Put $\A=\mathbb Z[v,v^{-1}]$ and $\mathbb K=\mathbb Q(v)$.
Subsections~\ref{subsec:reduced-graphs}--\ref{subsec:ascent-characteristics}
describe the two reduced Gelfand graphs and their insertion labels.
Subsections~\ref{subsec:quantum-based-prelim} and~\ref{subsec:levi-tops-prelim}
fix the based-module conventions used in
Section~\ref{sec:levi-top-criterion}.

\subsection{The reduced Gelfand graphs}
\label{subsec:reduced-graphs}

The type-$A_{d-1}$ Hecke algebra $\Hh_d$ is generated by
$H_1,\ldots,H_{d-1}$, subject to the braid relations and
\[
 (H_i-v)(H_i+v^{-1})=0.
\]
We use the standard Coxeter and Hecke-algebra conventions of
\cite{BjornerBrenti2005,KazhdanLusztig1979}.
Write $[d]=\{1,\ldots,d\}$, let $\I_d$ be the involutions in $S_d$,
and put
\[
 \Fix(w)=\{i\in[d]:w(i)=i\},\qquad
 \I_d^{\FPF}=\{w\in\I_d:\Fix(w)=\varnothing\}.
\]
For $w\in\I_d$ with fixed points $c_1<\cdots<c_k$, define two
fixed-point-free completions in $S_{2d}$.  Keep every two-cycle of $w$,
pair $c_j$ with $d+j$ for the ascending completion and with
$d+k+1-j$ for the descending completion, and pair the unused auxiliary
letters $d+k+1,\ldots,2d$ consecutively.  Denote the resulting involutions
by $\iota_{\mathrm{asc}}(w)$ and $\iota_{\mathrm{des}}(w)$.
The row and column standard modules are the free $\A$-modules with
bases $\{M_w:w\in\I_d\}$ and $\{N_w:w\in\I_d\}$, carrying the
ascending and descending completion labels, respectively.

Let $s_i=(i,i+1)$ and put $w'=s_iws_i$.  When $w'\ne w$, say that
$w'$ is above $w$ if
\[
 \ell(\iota_{\mathrm{asc}}(w'))>
 \ell(\iota_{\mathrm{asc}}(w));
\]
this is a strict ascent of $w$, and the reverse comparison is a strict
descent.  On the block with $k$ fixed points,
$\ell(\iota_{\mathrm{des}}(w))-
\ell(\iota_{\mathrm{asc}}(w))=k(k-1)$
\cite[Proposition~3.4]{MarbergZhangInsertion2024}, so the descending
completion gives the same strict comparisons.  The standard actions are
\begin{align}
 H_iM_w&=
 \begin{cases}
 M_{w'},&w'\text{ is above }w,\\
 M_{w'}+(v-v^{-1})M_w,&w'\text{ is below }w,\\
 -v^{-1}M_w,&w(i)=i\text{ and }w(i+1)=i+1,\\
 vM_w,&w(i)=i+1,
 \end{cases}
 \label{eq:four-cases-row}\\[2mm]
 H_iN_w&=
 \begin{cases}
 N_{w'},&w'\text{ is above }w,\\
 N_{w'}+(v-v^{-1})N_w,&w'\text{ is below }w,\\
 vN_w,&w(i)=i\text{ and }w(i+1)=i+1,\\
 -v^{-1}N_w,&w(i)=i+1.
 \end{cases}
 \label{eq:four-cases}
\end{align}
These actions and their compatible bar involutions are those of
\cite[Theorems~3.5 and~3.6 and Remark~3.7]{MarbergZhangInsertion2024};
see also \cite[Theorems~1.7 and~1.8]{MarbergZhangGelfand2022}.

Conjugation preserves the number of fixed points, so both standard
modules split into fixed-point blocks.  In the block with $r=2p$
moving points and $k=d-r$ fixed points, put
\[
 w_{r,k}=(1,2)(3,4)\cdots(r-1,r)(r+1)\cdots(d).
\]
Define the relative completion height and its strict order by
\begin{equation}\label{eq:height-order}
 h(w)=\frac{\ell(\iota_{\mathrm{asc}}(w))-
                    \ell(\iota_{\mathrm{asc}}(w_{r,k}))}{2},
 \qquad
 y\prec w\ \Longleftrightarrow\
 \#\Fix(y)=\#\Fix(w)\ \text{and}\ h(y)<h(w).
\end{equation}
A strict ascent increases the completion length by two and $h$ by one;
the descending completion gives the same height.

Every vertex is reached from $w_{r,k}$ by strict ascents.  Indeed, for
two adjacent positions whose conjugation is strict, a moving point
followed by a fixed point is an ascent, the reverse is a descent, and
two moving points give an ascent precisely when $w(i)<w(i+1)$.
These assertions follow by comparing their images in the ascending
completion.  A vertex with no strict descent therefore has all fixed
points at the right.  If its first remaining moving point $a$ is paired
with $b>a+1$, then $w(b-1)>w(b)=a$ gives a strict descent.  Removing
consecutive pairs inductively proves that $w_{r,k}$ is the unique vertex
with no strict descent.  Descending from any vertex and reversing the
path proves the assertion.  In particular, $h$ is a nonnegative integer.

On each block the compatible bar is characterized by
\[
 \overline v=v^{-1},\qquad
 \overline{H_ix}=H_i^{-1}\bar x,
\]
and by fixing $M_{w_{r,k}}$ or $N_{w_{r,k}}$.  Compatibility determines
the bar uniquely along the strict-ascent paths from this minimum.
It is triangular for $\prec$: write $q=v-v^{-1}$ and choose a strict
ascent $w=s_ius_i$.  Inductively,
\[
 \overline{X_w}=(H_i-q)\overline{X_u}
 \in X_w+\sum_{y\prec w}\A X_y,\qquad X=M,N.
\]
The leading term contributes $X_w-qX_u$; each other term starts at
height at most $h(u)-1$, and $H_i$ raises height by at most one.
The induction starts at the bar-fixed minimum.  Triangular induction
therefore gives the unique bar-invariant lower canonical bases
\begin{equation}\label{eq:canonical-normalization}
 \underline X_w\in
 X_w+\sum_{y\prec w}v^{-1}\mathbb Z[v^{-1}]X_y,
 \qquad X=M,N.
\end{equation}
Both bars and these canonical bases respect the fixed-point blocks.

The ascent labels $\mAsc(w)$ and $\nAsc(w)$ are read from the diagonal
$-v^{-1}$ coefficient in the respective canonical actions.
If $w'=s_iws_i\ne w$, then $i$ belongs to both labels precisely when
$w'$ is above $w$; in the weak cases,
\[
 \begin{array}{c|cc}
  &w(i)=i,\ w(i+1)=i+1&w(i)=i+1\\ \hline
  i\in\mAsc(w)&\text{yes}&\text{no}\\
  i\in\nAsc(w)&\text{no}&\text{yes}.
 \end{array}
\]
These are the labels in
\cite[equations~(3.6)--(3.7)]{MarbergZhangInsertion2024}.
Write $b_y$ for either canonical basis and $\Asc(y)$ for its label.
If $b_z=\sum_y p_{yz}X_y$ in the corresponding standard basis, put
\[
 \omega(y,z)=[v^{-1}](p_{yz}+p_{zy}).
\]
The canonical action is
\begin{equation}\label{eq:canonical-wgraph-action}
 H_i b_y=
 \begin{cases}
  vb_y,&i\notin\Asc(y),\\
  -v^{-1}b_y+\displaystyle\sum_{\substack{z\ne y\\i\notin\Asc(z)}}
       \omega(y,z)b_z,&i\in\Asc(y).
 \end{cases}
\end{equation}

The reduced row and column Gelfand graphs, denoted by $\mGamma$ and
$\nGamma$, have coefficient arrows: for distinct $y,z\in\I_d$, write
$y\longrightarrow_{\mcell}z$ when $\underline M_z$ has a nonzero
coefficient in some $H_i\underline M_y$, and write
$y\longrightarrow_{\ncell}z$ using $\underline N$ in the same way.
Equivalently, in either graph,
\begin{equation}\label{eq:reduced-arrow}
 y\longrightarrow z
 \quad\Longleftrightarrow\quad
 \omega(y,z)\ne0\quad\text{and}\quad
 \Asc(y)\nsubseteq\Asc(z).
\end{equation}
When vertices carry their completion labels, these arrows are written
$\iota_{\mathrm{asc}}(y)\longrightarrow_{\mcell}\iota_{\mathrm{asc}}(z)$
and
$\iota_{\mathrm{des}}(y)\longrightarrow_{\ncell}\iota_{\mathrm{des}}(z)$.
This is the reduced convention of
\cite[Theorem~1.10]{MarbergZhangGelfand2022}, with the column vertices
relabelled by descending completions.  The insertion paper retains the
symmetric weights $\omega$
\cite[Remark~3.1 and equation~(3.8)]{MarbergZhangInsertion2024};
the label test in \eqref{eq:canonical-wgraph-action} gives the same Hecke
action and the coefficient arrows \eqref{eq:reduced-arrow}.

The strongly connected components are the $\mcell$-cells and
$\ncell$-cells, respectively, and the connected components after retaining
only bidirected edges are the $\mcell$-molecules and $\ncell$-molecules.
No coefficient arrow crosses fixed-point blocks.  All cell and dominance
statements below use these reduced graphs.

\begin{remark}[Diagonal coefficients and reduced arrows]
\label{rem:reduced-coefficient-example}
Consider the column block with $d=3$ and one fixed point.  Put
$A=N_{(12)}$, $B=N_{(13)}$, and $C=N_{(23)}$.
The standard action and the canonical normalization
\eqref{eq:canonical-normalization} give
\[
 \underline N_{(12)}=A,\qquad
 e:=\underline N_{(13)}=B+v^{-1}A,\qquad
 f:=\underline N_{(23)}=C+v^{-1}B.
\]
Direct substitution into \eqref{eq:four-cases} gives
\[
 H_1e=-v^{-1}e+f,\quad H_1f=vf,\qquad
 H_2e=ve,\quad H_2f=e-v^{-1}f.
\]
The label of $e$ contains $1$, although $e$ is not an $H_1$ eigenvector.
This example also distinguishes the two weight conventions.  The
unreduced weights satisfy
$\omega((12),(13))=\omega((13),(12))=1$, while
\[
 \nAsc((12))=\{1,2\},\qquad \nAsc((13))=\{1\}.
\]
Indeed, $H_2A=-v^{-1}A+e$, whereas neither generator sends $e$ to a
vector with a nonzero $A$ coefficient.  Thus
$(12)\longrightarrow_{\ncell}(13)$ is a reduced arrow and its reverse
is absent.
\end{remark}

\subsection{Modified insertion and molecules}
\label{subsec:insertion-molecules}

A \emph{partially standard tableau} is a Young tableau with distinct
positive-integer entries increasing along rows and columns.
For a partition $\lambda\vdash N$,
write $\SYT(\lambda)$ for the standard Young tableaux of shape $\lambda$.
Ordinary Schensted row insertion of a new entry $a$ into $T$ is denoted
by $T\RSinsert a$.

\begin{definition}[Row and column Beissinger pair insertion]
Let $T$ be partially standard, and let $a\leq b$ be integers not occurring
in $T$.  If $a<b$ and $(i,j)$ is the new box of $T\RSinsert a$, define
$T\rBinsert(a,b)$ by appending $b$ to row $i+1$ of $T\RSinsert a$,
and define $T\cBinsert(a,b)$ by appending $b$ to the bottom of
column $j+1$ of the same tableau.
If $a=b$, append $a$ to the first row for $T\rBinsert(a,a)$ and to the
first column for $T\cBinsert(a,a)$.
\end{definition}

These are the operations of \cite[Algorithm~3.1]{Beissinger1987} and
\cite[Definition~2.10]{MarbergZhangInsertion2024}.
An individual operation need not preserve standardness for arbitrary
pairs; in the correspondences below, pairs are inserted in increasing
order of their second coordinates.

\begin{definition}[Full and restricted Beissinger correspondences]
For an involution $z$, let $P_{\mathrm{rB}}(z)$ and
$P_{\mathrm{cB}}(z)$ be obtained from the empty tableau by successive
row and column pair insertion, respectively, of all pairs $(a,b)$ with
$a\leq b=z(a)$, ordered by increasing $b$.
For $w\in I_d$, list its nontrivial cycles and fixed points as
\[
 (a_1,b_1),\ldots,(a_p,b_p),\qquad
 a_j<b_j,\quad b_1<\cdots<b_p,
 \qquad c_1<\cdots<c_k.
\]
Put
\begin{align*}
 U_{\mathrm r}(w)&=
 \emptyset\rBinsert(a_1,b_1)\cdots\rBinsert(a_p,b_p),&
 \PrB(w)&=U_{\mathrm r}(w)\RSinsert c_1\cdots\RSinsert c_k,\\
 U_{\mathrm c}(w)&=
 \emptyset\cBinsert(a_1,b_1)\cdots\cBinsert(a_p,b_p),&
 \PcB(w)&=U_{\mathrm c}(w)\RSinsert c_k\cdots\RSinsert c_1.
\end{align*}
The maps $\PrB$ and $\PcB$ are the restricted row and column
correspondences, and their shapes are denoted by
\[
 \lambda_{\mathrm r}(w)=\shape\PrB(w),\qquad
 \lambda_{\mathrm c}(w)=\shape\PcB(w).
\]
\end{definition}

For the two involutions in Remark~\ref{rem:reduced-coefficient-example},
the restricted column algorithm gives
\[
 \PcB((12))=\rowtableau{1&2&3\cr},\qquad
 \PcB((13))=\rowtableau{1&2\cr3\cr}.
\]
Their shapes are $(3)$ and $(2,1)$.  The reduced arrow goes from the
first shape to the second, in agreement with
Theorem~\ref{thm:main-intro}.

The following result relates these algorithms to the completion labels
of the two graphs.  For a tableau $T$ on an alphabet containing $[d]$,
write $T|_{[d]}$ for the tableau obtained by deleting entries greater
than $d$.

\begin{lemma}[Restricted Beissinger correspondences]
\label{lem:restricted-completion-bijection}
For every \(w\in I_d\),
\[
 \PrB(w)=P_{\mathrm{rB}}(\iota_{\mathrm{asc}}(w))|_{[d]},
 \qquad
 \PcB(w)=P_{\mathrm{cB}}(\iota_{\mathrm{des}}(w))|_{[d]}.
\]
For each \(k\), the first map is a bijection from the involutions with
\(k\) fixed points to the standard tableaux with \(k\) odd columns, while
the second is a bijection from the same set to the standard tableaux with
\(k\) odd rows.
\end{lemma}
\begin{proof}
The completion identities are
\cite[Lemmas~3.13 and~3.16]{MarbergZhangInsertion2024}, with the restricted
algorithms written in the equivalent form used above.  The bijectivity and
parity assertions are
\cite[Theorems~3.14 and~3.17]{MarbergZhangInsertion2024}.
\end{proof}

To describe the molecules, take $T\in\SYT(\lambda)$, where
$\lambda\vdash N$, and let $\operatorname{rr}(T)$ be its row-reading
word: read each row from left to right, starting with the bottom row.
Write $s_jT$ for the filling obtained by
interchanging $j$ and $j+1$.  For $2\leq i\leq N-1$, define
the elementary dual-equivalence involution
\begin{equation}\label{eq:dual-equivalence}
 D_i(T)=
 \begin{cases}
 s_{i-1}T,&i+1\text{ lies between }i\text{ and }i-1
                 \text{ in }\operatorname{rr}(T),\\
 s_iT,&i-1\text{ lies between }i\text{ and }i+1
                 \text{ in }\operatorname{rr}(T),\\
 T,&i\text{ lies between }i-1\text{ and }i+1
                 \text{ in }\operatorname{rr}(T).
 \end{cases}
\end{equation}
In either nontrivial case the exchanged boxes are incomparable in the Young
poset, so the new filling is standard; the same betweenness check gives
$D_i^2=1$.
These are Haiman's dual-equivalence moves \cite{Haiman1992}; their
graph-theoretic formulation is also standard in the theory of
dual-equivalence graphs \cite{Assaf2015}.

\begin{lemma}[Connectivity of elementary dual equivalence]
\label{lem:dual-equivalence-connected}
For every partition $\lambda\vdash N$, the graph on $\SYT(\lambda)$ joining
$T$ to $D_i(T)\neq T$ for all $2\leq i\leq N-1$ is connected.
\end{lemma}
\begin{proof}
This is the standard consequence of the Knuth and dual-Knuth
characterization in \cite[Theorem~2.3]{MarbergZhangInsertion2024}, explained
immediately after that theorem; see also \cite{Haiman1992}.
\end{proof}

\begin{proposition}[Molecules in the reduced Gelfand graphs]
\label{prop:reduced-molecules}
The bidirected coefficient edges of $\mGamma$ and $\nGamma$ correspond,
respectively, to elementary dual-equivalence moves on $\PrB$ and $\PcB$.
Consequently the molecules of each reduced graph are exactly the fibers
of its modified insertion shape.
\end{proposition}
\begin{proof}
For the reduced graphs, the bidirected edges are characterized by
\cite[Corollaries~3.14 and~3.17 and Lemma~3.27]{MarbergZhangGelfand2022}.
The local reduction in the proof of
\cite[Lemma~3.10]{MarbergZhangInsertion2024} expresses these criteria
as the row and column conjugation relations used there.
The local equivalences in
\cite[Theorems~3.15 and~3.18]{MarbergZhangInsertion2024}
identify those relations with the moves $D_i$ on the restricted tableaux.
Lemma~\ref{lem:restricted-completion-bijection} and
Lemma~\ref{lem:dual-equivalence-connected} now identify their connected
components with shape fibers.
\end{proof}

\subsection{Ascent labels and Schur characteristics}
\label{subsec:ascent-characteristics}

The ascent labels encode the Schur characteristic of each insertion-shape
fiber.  For $T\in\SYT(\lambda)$ with $\lambda\vdash N$, put
\[
 \Des(T)=\{i\in[N-1]:i+1\text{ lies in a lower row than }i\}.
\]
Write $s_\lambda$ for the Schur function.  For $D\subseteq[N-1]$,
Gessel's fundamental quasisymmetric function is
\[
 F_D=\sum_{\substack{1\leq a_1\leq\cdots\leq a_N\\
                     j\in D\Rightarrow a_j<a_{j+1}}}
       x_{a_1}\cdots x_{a_N}.
\]
For every partition $\mu\vdash N$, the elementary standardization
argument for semistandard tableaux gives Gessel's expansion
\cite{Gessel1984}:
\begin{equation}\label{eq:schur-fundamental}
 s_\mu=\sum_{S\in\SYT(\mu)}F_{\Des(S)}.
\end{equation}
Indeed, standardize a semistandard tableau by reading equal entries from left
to right.  Its weakly increasing content word has forced strict rises exactly
at the descents of the resulting standard tableau.  This partitions all
semistandard tableaux and proves \eqref{eq:schur-fundamental} coefficient by
coefficient.

\begin{proposition}[Ascent labels and insertion-shape characteristics]
\label{prop:insertion-shape-characteristics}
\label{prop:ascent-transposed-descent}
\label{lem:insertion-carrier-graphs}
For every \(w\in I_d\),
\[
 \mAsc(w)=\Des(\PrB(w)^{\mathsf t}),
 \qquad
 \nAsc(w)=\Des(\PcB(w)^{\mathsf t}).
\]
For each admissible fixed-point block and each shape $\lambda$
occurring in that block,
\begin{align}
 \sum_{\lambda_{\mathrm r}(w)=\lambda}F_{\mAsc(w)}
   &=s_{\lambda^{\tr}},\label{eq:row-shape-character-new}\\
 \sum_{\lambda_{\mathrm c}(w)=\lambda}F_{\nAsc(w)}
   &=s_{\lambda^{\tr}}.\label{eq:column-shape-character-new}
\end{align}
Each sum is over the fixed-point block under consideration.
\end{proposition}
\begin{proof}
The ascent formulas and adjacent-letter calculations of
\cite[equations~(3.6)--(3.7), Proposition~3.4, and the proofs of
Theorems~3.15 and~3.18]{MarbergZhangInsertion2024} give
\[
 \mAsc(w)=[d-1]\setminus\Des(\PrB(w)),
 \qquad
 \nAsc(w)=[d-1]\setminus\Des(\PcB(w)).
\]
For every standard tableau \(T\), exactly one of \(i+1\) lying below \(i\)
in \(T\) and in \(T^{\mathsf t}\) occurs.  Hence
\(\Des(T^{\mathsf t})=[d-1]\setminus\Des(T)\), proving the ascent
identities.

Lemma~\ref{lem:restricted-completion-bijection} identifies each shape
fiber with all standard tableaux of that shape.  The ascent identities
therefore identify its labels with the descents of the transposed tableaux.
Applying \eqref{eq:schur-fundamental} to $\lambda^{\tr}$ gives
\[
 s_{\lambda^{\tr}}=
 \sum_{T\in\SYT(\lambda)}F_{\Des(T^{\tr})},
\]
which proves both characteristic identities.
\end{proof}

\subsection{Quantum groups and lower global bases}
\label{subsec:quantum-based-prelim}

Sections~\ref{sec:row-model} and~\ref{sec:column-model} compare the Schur
characteristics of Proposition~\ref{prop:insertion-shape-characteristics}
with highest-weight labels.  This comparison uses realizations of the
two canonical Hecke bases inside lower based modules.  We fix the
required quantum-group and global-basis conventions here.

Fix a finite-type root datum with simple roots $\alpha_i$, simple
coroots $h_i$, and Cartan matrix $a_{ij}=\langle h_i,\alpha_j\rangle$
indexed by $I$.  A weight $\lambda$ is \emph{dominant} if
$\langle h_i,\lambda\rangle\geq0$ for all $i\in I$.
Choose positive integers $d_i$ with $d_i a_{ij}=d_j a_{ji}$ and put
$v_i=v^{d_i}$.  The quantum integers and factorials are
\[
 [a]_i=\frac{v_i^a-v_i^{-a}}{v_i-v_i^{-1}},\qquad
 [a]_i!=\prod_{r=1}^a[r]_i,\qquad [0]_i!=1.
\]
The Drinfeld--Jimbo algebra $U_v(\mathfrak g)$ is generated over
$\mathbb K$ by $E_i,F_i,K_i^{\pm1}$, with relations
\begin{gather*}
 K_iK_i^{-1}=K_i^{-1}K_i=1,\qquad K_iK_j=K_jK_i,\\
 K_iE_jK_i^{-1}=v_i^{a_{ij}}E_j,\qquad
 K_iF_jK_i^{-1}=v_i^{-a_{ij}}F_j,\\
 E_iF_j-F_jE_i=\delta_{ij}\frac{K_i-K_i^{-1}}{v_i-v_i^{-1}}.
\end{gather*}
Writing
\[
 E_i^{(a)}=\frac{E_i^a}{[a]_i!},\qquad
 F_i^{(a)}=\frac{F_i^a}{[a]_i!}
\]
for the divided powers, the remaining defining relations are
\begin{equation}\label{eq:quantum-serre-relations}
 \sum_{r=0}^{1-a_{ij}}(-1)^r
 E_i^{(r)}E_jE_i^{(1-a_{ij}-r)}=0,
 \qquad
 \sum_{r=0}^{1-a_{ij}}(-1)^r
 F_i^{(r)}F_jF_i^{(1-a_{ij}-r)}=0
 \qquad(i\ne j).
\end{equation}
Tensor products use the coproduct
\begin{equation}\label{eq:coproduct-convention}
 \begin{aligned}
 \Delta(K_i)&=K_i\otimes K_i,\\
 \Delta(E_i)&=E_i\otimes K_i+1\otimes E_i,\\
 \Delta(F_i)&=F_i\otimes 1+K_i^{-1}\otimes F_i.
 \end{aligned}
\end{equation}

All based modules below are type-$1$ integrable with finite-dimensional
weight spaces: they decompose into the spaces
\[
 X_\lambda=\{x\in X:K_ix=v_i^{\langle h_i,\lambda\rangle}x
 \text{ for all }i\in I\},
\]
and $E_i,F_i$ act locally nilpotently.  The $i$-string decomposition
writes each weight vector uniquely as
$x=\sum_{a\geq0}F_i^{(a)}x_a$, where $E_ix_a=0$.
Kashiwara's operators are
\[
 \widetilde e_i x=\sum_{a\geq1}F_i^{(a-1)}x_a,
 \qquad
 \widetilde f_i x=\sum_{a\geq0}F_i^{(a+1)}x_a.
\]

Write $U_\A(\mathfrak g)$ for Lusztig's divided-power integral form.
An integral lattice $X_\A$ is a free $U_\A$-stable $\A$-submodule
with $\mathbb K\otimes_\A X_\A=X$.  Put
\[
 \Oinf=\mathbb Z[v^{-1}]_{\{g:g(0)=1\}},
\]
the localization obtained by inverting polynomials with constant
term $1$.  A lower crystal lattice $\mathcal L^-$ is a free
$\Oinf$-lattice stable under $\widetilde e_i,\widetilde f_i$.
A crystal basis $B(X)$ is a weight-compatible $\mathbb Z$-basis of
$\mathcal L^-/v^{-1}\mathcal L^-$ on which these operators preserve
$B(X)\sqcup\{0\}$ and are partial inverses.  Write
$b\xrightarrow{i}b'$ when $\widetilde f_i b=b'$, and set
\[
 \varepsilon_i(b)=\max\{a:\widetilde e_i^ab\ne0\},\qquad
 \varphi_i(b)=\max\{a:\widetilde f_i^ab\ne0\}.
\]
We use normal highest-weight crystals: restriction to any finite-type
set of colors is a disjoint union of crystals of integrable
highest-weight modules.  Each such component has a unique highest
element.

The algebra and each based module carry a semilinear bar involution
with $\overline v=v^{-1}$,
$\overline{E_i}=E_i$, $\overline{F_i}=F_i$, and
$\overline{K_i}=K_i^{-1}$.  Put $\mathcal L^+=\overline{\mathcal L^-}$.
Put $\mathcal E=X_\A\cap\mathcal L^-\cap\mathcal L^+$.
Our integral \emph{balanced} structure requires that reduction
\[
 \mathcal E\longrightarrow\mathcal L^-/v^{-1}\mathcal L^-
\]
be an isomorphism of $\mathbb Z$-modules.  We also require the
inverse images $G(b)$ of $b\in B(X)$ to give the three lattice bases
\[
 X_\A=\bigoplus_{b\in B(X)}\A G(b),\qquad
 \mathcal L^-=\bigoplus_{b\in B(X)}\Oinf G(b),\qquad
 \mathcal L^+=\bigoplus_{b\in B(X)}\overline{\Oinf}\,G(b).
\]
We assume that the lattices respect weights and that every $G(b)$
is bar invariant.  The set
$\mathbf B(X)=\{G(b):b\in B(X)\}$ is the \emph{lower global basis},
and $(X,\mathbf B(X))$ is called a \emph{lower based module}
\cite{Kashiwara1991,Kashiwara1993,Lusztig1993}.
Balancedness gives the uniqueness relation
\begin{equation}\label{eq:balanced-uniqueness-prelim}
 X_\A\cap v^{-1}\mathcal L^-\cap v\mathcal L^+=0.
\end{equation}

Write $U_v^-$ for the subalgebra generated by the $F_i$.  For $i\in I$
and $a\geq1$, the ideal $U_v^-F_i^{(a)}$ is spanned over $\mathbb K$
by the lower global basis vectors $G(b)$ for which $\widetilde e_i^a$
is nonzero on the image of $b$ under Kashiwara's involution
\cite[Theorem~2.29]{Kimura2012}.  This involution on crystal labels
is induced by the antiautomorphism that fixes each $F_i$ and reverses
products.  For a dominant weight $\lambda$, let $u_\lambda$ be the
highest global basis vector of $V(\lambda)$.  Under the projection
$x\mapsto xu_\lambda$, each $G(b)$ maps to zero or to the lower global
basis vector indexed by the image of $b$ in the quotient crystal,
with coefficient one \cite[Lemma~7.3.2]{Kashiwara1991}
(see also \cite[Theorem~2.26]{Kimura2012}).

\subsection{Levi tops and dominance}
\label{subsec:levi-tops-prelim}
\label{subsec:based-levi-prelim}

Restricting a lower based module to a Levi subalgebra supplies the
highest-weight labels used in Theorem~\ref{thm:levi-top-criterion}.
For $J\subseteq I$, let $U_v(\mathfrak g_J)$ be the derived Levi
quantum subalgebra generated by $E_j,F_j,K_j^{\pm1}$ for $j\in J$.
We record the ambient weight grading separately.  In type $A$ we
recover the polynomial weight from the derived Levi weight and the
total degree, as in Subsection~\ref{subsec:squarefree-tools}.
The $J$-crystal retains only the arrows colored by $J$.
By normality, the component containing $b$ has a unique highest
element; its dominant $J$-weight is the \emph{lower Levi top}
$\topJ(b)$.
Write $V_J(\lambda)$ for the simple integrable Levi module of
highest weight $\lambda$ and $B_J(\lambda)$ for its crystal.
The restriction is \emph{multiplicity-free} if each highest-weight
crystal $B_J(\lambda)$ occurs at most once.

For dominant weights write $\alpha\preceq\beta$ if $\beta-\alpha$
is a nonnegative integral combination of simple roots.
For partitions of the same integer, this gives the usual dominance order,
\[
 \alpha\preceq\beta
 \quad\Longleftrightarrow\quad
 \alpha_1+\cdots+\alpha_j\leq\beta_1+\cdots+\beta_j
 \qquad(j\geq1).
\]
Transposition reverses this order:
\begin{equation}\label{eq:transpose-reverses}
 \alpha\preceq\beta\quad\Longleftrightarrow\quad
 \beta^{\tr}\preceq\alpha^{\tr}.
\end{equation}

An \emph{upset} of a poset is closed upward, and a \emph{downset}
is closed downward.  A submodule $L'$ of a free $\A$-module $L$
is \emph{saturated} if $L/L'$ is torsion-free.

Let $\Lambda$ be the set of tops occurring in a lower based module $M$.
For $\Omega\subseteq\Lambda$ write
\[
 M_\Omega=\bigoplus_{\topJ(b)\in\Omega}\A G(b).
\]
For $\mu\in\Lambda$ let $\Omega_\mu=\{\nu\in\Lambda:\mu\preceq\nu\}$.
The \emph{one-top layer} at $\mu$ is
$M_{\Omega_\mu}/M_{\Omega_\mu\setminus\{\mu\}}$.
Its basis consists of the images of the vectors of top $\mu$.
Theorem~\ref{thm:excellent} proves that both spans are Levi submodules
and that the quotient contains one based copy of $V_J(\mu)$ for
each crystal component with that top.

In the next section, Theorem~\ref{thm:levi-top-criterion} identifies
cells from multiplicity-free Levi tops, stability of upset spans under
the coefficient operators, and bidirected connectivity of the fibers.
Theorem~\ref{thm:excellent} supplies the Levi stability used in the
two applications.

\section{General criteria and based-module tools}
\label{sec:levi-top-criterion}

We collect the graph-theoretic criterion and the based-module tools
used for both models.

\subsection{The graph-theoretic criterion}

Let $\mathscr H$ be an $\A$-algebra, and fix elements
$h_s\in\mathscr H$ for $s\in\mathscr S$.  Let $\mathscr M$ be an
$\mathscr H$-module whose underlying $\A$-module is free with basis
$\{C_b:b\in\mathscr B\}$.  For distinct
$b,c\in\mathscr B$, write $b\longrightarrow c$ when $C_c$ has nonzero
coefficient in $h_sC_b$ for some $s\in\mathscr S$.  A \emph{cell} is a
strongly connected component of this directed graph, and a \emph{molecule}
is a connected component after only the bidirected edges
$b\longleftrightarrow c$ are retained.

Suppose the basis comes from a based Levi restriction, or a specified
weight slice of such a restriction, as in
Subsection~\ref{subsec:levi-tops-prelim}.  A
\emph{multiplicity-free Levi-top datum} consists of a finite poset
$(\Lambda,\preceq)$ of distinct dominant Levi weights and a surjective map
\[
 \tau:\mathscr B\longrightarrow\Lambda
\]
such that the basis vectors with top $\lambda$ correspond to the
global basis of the unique one-top layer isomorphic to $V_J(\lambda)$,
or to the global basis of its specified weight slice.
In the row and column applications, $\preceq$ is the dominance order
of Subsection~\ref{subsec:levi-tops-prelim}.
Thus the top distinguishes the unique highest-weight constituent of
that layer, and no additional label for a multiplicity space is needed.
Write $\mathscr B_\lambda=\tau^{-1}(\lambda)$, and, for
$\Omega\subseteq\Lambda$, put
\begin{equation}\label{eq:abstract-top-span}
 \mathscr M_\Omega=
 \bigoplus_{\tau(b)\in\Omega}\A C_b.
\end{equation}

\begin{theorem}[Multiplicity-free Levi-top criterion]
\label{thm:levi-top-criterion}
Assume that $\tau:\mathscr B\to\Lambda$ is a multiplicity-free Levi-top
datum and that $\mathscr M_\Omega$ is an $\mathscr H$-submodule
for every upset $\Omega\subseteq\Lambda$.
If, in addition, the graph induced on every fiber
$\mathscr B_\lambda$ is connected by bidirected edges, then
\[
 \boxed{\text{cells}=\text{molecules}=\text{Levi-top fibers}.}
\]
More precisely, the three collections of subsets of $\mathscr B$ are all
$\{\mathscr B_\lambda:\lambda\in\Lambda\}$.
\end{theorem}

\begin{proof}
If $b\longrightarrow c$,
take the principal upset
\(
 \Omega_b=\{\lambda\in\Lambda:\tau(b)\preceq\lambda\}.
\)
Since $C_b\in\mathscr M_{\Omega_b}$ and this span is a submodule, every
nonzero coefficient in every $h_sC_b$ is still indexed by
$\Omega_b$.  Hence
\begin{equation}\label{eq:top-monotone-up}
 b\longrightarrow c\quad\Longrightarrow\quad
 \tau(b)\preceq\tau(c).
\end{equation}
For a bidirected edge, applying \eqref{eq:top-monotone-up} in both directions
and using antisymmetry gives $\tau(b)=\tau(c)$.  Thus every molecule is
contained in a single Levi-top fiber.

If $b$ and $c$ lie in the same cell, there are directed paths from $b$ to
$c$ and from $c$ to $b$.  Iterating \eqref{eq:top-monotone-up} along both
paths gives
$\tau(b)\preceq\tau(c)\preceq\tau(b)$, so again
$\tau(b)=\tau(c)$.  Therefore every cell is contained in one fiber.
Conversely, the connectivity hypothesis puts the whole fiber
$\mathscr B_\lambda$ in one molecule.  Since all bidirected edges preserve
$\tau$, that molecule contains nothing outside the fiber.  Hence each fiber
is exactly one molecule.  A molecule is contained in a cell, while every
cell is contained in a fiber, so all three subsets coincide.

\end{proof}

To make use of Theorem~\ref{thm:levi-top-criterion}, we define a model and verify the following conditions.
\begin{enumerate}[label=\textup{(C\arabic*)},ref=C\arabic*,noitemsep]
\item \label{cond:levi-top-multiplicity}
The map $\tau:\mathscr B\to\Lambda$ is a multiplicity-free Levi-top datum.
\item \label{cond:levi-top-upsets}
The span $\mathscr M_\Omega$ is an $\mathscr H$-submodule
for every upset $\Omega\subseteq\Lambda$.
\item \label{cond:levi-top-connectivity}
Every fiber $\mathscr B_\lambda$ is connected by bidirected coefficient arrows.
\end{enumerate}

\begin{remark}
Multiplicity-free branching does not imply
conditions~\ref{cond:levi-top-upsets} and
\ref{cond:levi-top-connectivity}.
With multiplicities, the top remains well defined but may label
several Levi components; distinguishing these requires additional
labels compatible with the based structure and the operators $h_s$,
together with separate verification of stability and connectivity.
\end{remark}

\subsection{Excellent restriction as a tool for stability}
\label{subsec:excellent-restriction}

The following integral form of excellent restriction
\cite{Kashiwara1991,Kashiwara1993,Lusztig1993} supplies the Levi
stability used to verify condition~\ref{cond:levi-top-upsets}.
For an $\A$-submodule $K\subseteq L$ write
\[
 \operatorname{Sat}_{\A}(K)=
 \{x\in L:\text{there is }0\ne a\in\A\text{ with }ax\in K\}.
\]

\begin{theorem}[Excellent restriction]\label{thm:excellent}
Let $M$ be a finite-dimensional integrable lower based module and let $J$ be
of finite type.  If $\Omega$ is an upset in the dominance order on dominant
$J$-weights, then
\begin{equation}\label{eq:excellent-span}
 M_\Omega=\bigoplus_{\topJ(b)\in\Omega}\A G(b)
\end{equation}
is a saturated $U_\A(\mathfrak g_J)$-submodule of $M_\A$.  If two consecutive upsets
differ by one weight $\mu$, the corresponding subquotient is a direct sum of
based highest-weight modules $V_J(\mu)$.
\end{theorem}
\begin{proof}
Normality expresses the character of the restriction as the sum of
the highest-weight characters of its $J$-crystal components.  Since the
finite-dimensional restriction is semisimple, its simple constituents
therefore have exactly these highest weights, with the same multiplicities.
Choose one maximal occurring $J$-top $\mu$, and let $b^+$ be a
$J$-highest crystal element of that top.
Although $\widetilde e_i b^+=0$, a priori $E_iG(b^+)$ could contain terms
invisible in the crystal reduction.  Every such term has weight
$\mu+\alpha_i$.  In a finite-dimensional integrable module a vector of that
weight can occur only in a simple $J$-constituent of highest weight $\nu$
with
\[
 \nu-(\mu+\alpha_i)\in\sum_{j\in J}\mathbb Z_{\geq0}\alpha_j.
\]
Thus $\nu$ strictly dominates $\mu$, contrary to maximality.  Hence
$E_iG(b^+)=0$ for every $i\in J$: the global basis lift of a maximal crystal
highest element is an actual highest-weight vector.

There is consequently a homomorphism
\[
 V_J(\mu)\longrightarrow M,
 \qquad v_\mu\longmapsto G(b^+).
\]
Denote this map by $\iota$.  It is injective over $\mathbb K$, since
$V_J(\mu)$ is simple.  It preserves the integral lattice because
$V_J(\mu)_\A=U_\A^-(\mathfrak g_J)v_\mu$ and $G(b^+)\in M_\A$.
It also commutes with bar, as its highest vector has bar-invariant image.

Applying $\iota$ to the $j$-string decomposition of a vector gives
the $j$-string decomposition of its image.  Thus $\iota$ commutes
with $\widetilde e_j$ and $\widetilde f_j$ for every $j\in J$.
The highest-weight lower lattice is the $\Oinf$-span of the vectors
$\widetilde f_{j_t}\cdots\widetilde f_{j_1}v_\mu$.
Their images lie in $\mathcal L^-$, since $G(b^+)\in\mathcal L^-$ and
$\mathcal L^-$ is stable under the Kashiwara operators.  Thus $\iota$
preserves the lower lattice, and its reduction identifies $B_J(\mu)$
with the normal $J$-component of $b^+$.  Bar compatibility also gives
preservation of the conjugate lattice.  Consequently the image of each
global basis vector of $V_J(\mu)$ and the corresponding $G(b)$ are
integral, bar invariant, and have the same crystal residue.  Their
difference lies in
$M_\A\cap v^{-1}\mathcal L^-\cap v\mathcal L^+$
and vanishes by \eqref{eq:balanced-uniqueness-prelim}.
Hence $\iota$ is a based embedding, with precisely the global basis
vectors in that $J$-component as its image basis.

Apply this construction to all components of the chosen top $\mu$.
Their image bases are disjoint subsets of the global basis of $M$, so
their sum is a direct sum of based highest-weight submodules.  Its basis
is exactly the set of $b$ with $\topJ(b)=\mu$.
The quotient is again a lower based integrable
$U_v(\mathfrak g_J)$-module.  Its integral lattice and its two crystal
lattices are obtained by deleting this basis subset.  The quotient
map commutes with the Kashiwara operators by the same string-decomposition
argument, so its crystal is obtained by deleting these entire
$J$-components and remains normal.  Choose one maximal top in the
quotient and repeat.
This constructs a filtration in any descending linear extension of the
occurring tops, with one weight removed at each step.  Each successive
quotient is a direct sum of based highest-weight modules for that weight.

Finally, given an arbitrary dominance upset $\Omega$, choose a descending
linear extension of the finite poset of occurring tops in which all elements
of $\Omega$ precede its complement; this is possible precisely because
$\Omega$ is an upset.  Stop the preceding construction after the last member
of $\Omega$.  The resulting submodule is \eqref{eq:excellent-span}.
The integral basis condition in Subsection~\ref{subsec:quantum-based-prelim}
identifies the quotient lattice with the free $\A$-span of the
complementary vectors.  Thus $M_\Omega$ is saturated.  The induction also identifies each one-top
subquotient.
\end{proof}

Suppose a based comparison realizes a Hecke block on $(M_\A)_\delta$
with generators $\mathsf H_i=E_iF_i-v^{-1}\operatorname{id}$ for $i\in J$.
By Theorem~\ref{thm:excellent}, $M_\Omega$ is stable under $E_i,F_i$;
since $E_iF_i$ has weight zero,
\[
 (M_\Omega)_\delta
 =\bigoplus_{\substack{\wt(b)=\delta\\\topJ(b)\in\Omega}}\A G(b)
\]
is Hecke-stable.  Transporting this basis subset verifies
condition~\ref{cond:levi-top-upsets}; its free complementary basis
also proves saturation.

Based one-top subquotients leave coefficients between surviving
basis vectors unchanged.  Thus bidirected edges established in a
layer remain bidirected in the original graph, as required when
verifying condition~\ref{cond:levi-top-connectivity}.

\subsection{Based quotients and braid conventions}
\label{subsec:quotient-conventions}

Let $U_v(\mathfrak g)$ have finite-type root datum, and let $w_0$ be
the longest Weyl-group element.  We use the symmetrizing parameters
and divided powers of Subsection~\ref{subsec:quantum-based-prelim} and
write $U_\A^+$ for the divided-power integral form of the positive half.
For a dominant weight $\Lambda$, let $v_{w_0\Lambda}$ be the normalized
lowest global basis vector of $V(\Lambda)$.  Write
$\pi_\Lambda:U_v^+\longrightarrow V(\Lambda)$ for the action map and put
\[
 \pi_\Lambda(x)=xv_{w_0\Lambda},\qquad
 m_i=\langle h_i,-w_0\Lambda\rangle.
\]
Thus $m_i$ is the length of the $i$-string starting at the lowest vector.

We identify the crystal of $U_v^+$ with $B(\infty)$ of the negative
half by $q\mapsto v^{-1}$ and $f_i\mapsto E_i$.  Let $\widetilde e_i$
be the operator that removes an $i$ from a crystal string, and let $*$
be the crystal involution induced by the algebra antiautomorphism fixing
$v$ and the $E_i$.  Since $*$ reverses products, the string parameter
for divided powers on the right is
\[
 \varepsilon_i^*(b)=\max\{a\geq0:(\widetilde e_i^*)^a b\ne0\},
 \qquad \widetilde e_i^*=*\widetilde e_i*.
\]

\begin{lemma}[Based extremal quotient]\label{lem:extremal-quotient}
For every dominant weight $\Lambda$,
\begin{equation}\label{eq:extremal-kernel}
 \ker(\pi_\Lambda|_{U_\A^+})
 =\bigoplus_{\exists i:\,\varepsilon_i^*(b)>m_i}\A G(b).
\end{equation}
Every basis element outside this kernel maps, with coefficient one,
to the global basis element having the same quotient-crystal label.
\end{lemma}
\begin{proof}
The lowest-weight presentation gives the rational kernel
\[
 \ker\pi_\Lambda=\sum_i U_v^+E_i^{(m_i+1)}.
\]
Under $q\mapsto v^{-1}$ and $f_i\mapsto E_i$, the divided-power ideal
property recalled in Subsection~\ref{subsec:quantum-based-prelim} gives
\[
 U_v^+E_i^{(a)}
 =\bigoplus_{\varepsilon_i^*(b)\geq a}\mathbb K G(b)
 \qquad(a\geq1).
\]
Taking $a=m_i+1$ and summing over $i$ shows that the rational kernel
is spanned by those $G(b)$ with $\varepsilon_i^*(b)>m_i$ for some $i$.
Since the global basis is an $\A$-basis of $U_\A^+$, intersecting this
span with $U_\A^+$ gives exactly the direct sum in
\eqref{eq:extremal-kernel}.  The quotient is free on the remaining
basis classes, so the kernel is saturated.

For the image of a surviving $G(b)$, apply the highest-weight projection
property of Subsection~\ref{subsec:quantum-based-prelim} to
$V(-w_0\Lambda)$.  The Chevalley exchange takes its normalized highest
vector to $v_{w_0\Lambda}$ in $V(\Lambda)$, so $\pi_\Lambda(G(b))$ is
the global basis vector with the quotient-crystal label of $b$, with
coefficient one.
\end{proof}

We use Lusztig's braid operators with the convention
\begin{equation}\label{eq:shared-braid-convention}
 \mathsf T_i=T'_{i,-1},\qquad
 \mathsf T_i(E_i)=-K_i^{-1}F_i,\qquad
 \mathsf T_i(F_i)=-E_iK_i,
\end{equation}
and, for $j\ne i$,
\[
 \mathsf T_i(E_j)=
 \sum_{r+s=-a_{ij}}(-1)^r v_i^{-r}
       E_i^{(r)}E_jE_i^{(s)}.
\]
These conventions apply to both root systems below, with their
respective symmetrizing parameters.  For PBW root vectors, convex
orders, quantum Schubert cells, and their canonical bases, see
\cite{Kimura2012,LevendorskiiSoibelman1991,Lusztig1993}.

Fix a reduced expression $w_0=s_{i_1}\cdots s_{i_N}$ and put
$\beta_k=s_{i_1}\cdots s_{i_{k-1}}(\alpha_{i_k})$.
For $a\geq0$ write
$E_{\beta_k}^{(a)}=\mathsf T_{i_1}\cdots\mathsf T_{i_{k-1}}(E_{i_k}^{(a)})$.
The forward divided PBW product with exponent tuple
$\mathbf b=(b_1,\ldots,b_N)\in\mathbb Z_{\geq0}^N$ is
\[
 P(\mathbf b)=E_{\beta_1}^{(b_1)}\cdots E_{\beta_N}^{(b_N)}.
\]
We write $\mathbf b<_{\mathrm{lex}}\mathbf c$ if, at the first coordinate
$k$ where the tuples differ, $b_k<c_k$.

\begin{lemma}[The ordinary bar in a forward PBW basis]
\label{lem:full-pbw-bar}
For every exponent tuple $\mathbf b$, ordinary bar satisfies
\begin{equation}\label{eq:full-pbw-ordinary-bar}
 \overline{P(\mathbf b)}\in P(\mathbf b)+
 \sum_{\substack{\mathbf c>_{\mathrm{lex}}\mathbf b\\
                  \wt P(\mathbf c)=\wt P(\mathbf b)}}
       \A P(\mathbf c).
\end{equation}
\end{lemma}
\begin{proof}
Use first the negative half with parameter $q$.  Write $F_{\mathbf b}$
for its forward divided PBW vectors formed with $T''_{i,1}$.
On each weight space let $D_{\mathbf c}$ be the basis dual to these
vectors for Kashiwara's bilinear form $(\, ,\, )_K$.
These are Kimura's dual PBW vectors
\cite[\S4.7.2]{Kimura2012}.  By his definition of the dual bar
in \S3.1 and its triangular expansion in Proposition~4.25,
\[
 (\sigma x,y)_K=\overline{(x,\bar y)_K},\qquad
 \sigma D_{\mathbf c}
   =D_{\mathbf c}+\sum_{\mathbf a<_{\mathrm{lex}}\mathbf c}
                  u_{\mathbf a,\mathbf c}D_{\mathbf a},
\]
where the order is that of \cite[\S4.3.5]{Kimura2012}.
The coefficient of $F_{\mathbf c}$ in $\overline{F_{\mathbf b}}$
is therefore $\overline{(\sigma D_{\mathbf c},F_{\mathbf b})_K}$.
It is one for $\mathbf c=\mathbf b$ and can be nonzero off the
diagonal only for $\mathbf b<_{\mathrm{lex}}\mathbf c$.
All coefficients belong to $\mathbb Z[q,q^{-1}]$ because ordinary
bar preserves the integral form and the full divided PBW products
are its basis \cite[Chapters~40--41]{Lusztig1993}.

The substitution $q\mapsto v^{-1}$, $f_i\mapsto E_i$,
$e_i\mapsto F_i$, $t_i\mapsto K_i$ commutes with ordinary bar.
The generator formulas for $T''_{i,1}$ in \cite[\S4.2.1]{Kimura2012}
become \eqref{eq:shared-braid-convention} and the corresponding
formula for $\mathsf T_i(E_j)$.  Thus $F_{\mathbf b}$ maps to
$P(\mathbf b)$ with the factors in the same order.
Bar preserves root degree, giving \eqref{eq:full-pbw-ordinary-bar}.
\end{proof}

We also use the integral and crystal PBW statements for the full word.
The divided products form an $\A$-basis of the integral positive half.
In the negative half at $q=0$, their residues are the distinct crystal
labels \cite[Theorems~4.9 and~4.16]{Kimura2012}.  They exhaust each
weight space since both bases have the same finite cardinality there.
To see the lattice statement over our integer coefficient ring,
order the lower global basis by these residues.  In that basis the
PBW transition matrix has entries in $\mathbb Z[q,q^{-1}]$ by
integrality, is regular at zero, and has value $I$ there.  Its
entries thus lie in $\mathbb Z[q]$.  Its determinant is a unit
of $\mathbb Z[q,q^{-1}]$ with value one at zero, hence equals one.
The inverse matrix also has entries in $\mathbb Z[q]$.
After $q=v^{-1}$ the divided products therefore form an
$\Oinf$-basis of the lower lattice with the specified residues.

\subsection{Squarefree highest-weight slices}
\label{subsec:squarefree-tools}

For the type-$A$ constructions, let $d\geq2$ and put
$J=\{1,\ldots,d-1\}$.  A polynomial module of degree $m$ is assigned
its $\mathfrak{gl}_d$ weight lift of total degree $m$, so its dominant
weights are partitions of $m$ with at most $d$ parts.  For the ambient Levi of rank $R-1$ we use the same lift in
$R$ coordinates and write $\mathbf 1_R=\sum_{i=1}^R\varepsilon_i$.
Write $V_m$ for the vector representation of the type-$A_{m-1}$ Levi, with
weight basis $e_1,\ldots,e_m$ and $\wt(e_i)=\varepsilon_i$.
In particular, $V(k\omega_1)=\operatorname{Sym}_v^k(V_d)$.
The two constructions use the coproduct \eqref{eq:coproduct-convention}
and the lower crystal lattice at $v=\infty$.  We specify the tensor-product
bar in this convention.  Put $Q_+=\sum_i\mathbb Z_{\geq0}\alpha_i$ and let
\[
 \Theta=\sum_{\beta\in Q_+}\Theta_\beta,\qquad
 \Theta_0=1\otimes1,\qquad
 \Theta_\beta\in U_v^+{}_{\beta}\otimes U_v^-{}_{-\beta}
\]
be the quasi-$R$-matrix with the normalization
\begin{equation}\label{eq:quasi-r-convention}
 \Theta\,\overline{\Delta(u)}=\Delta(\bar u)\,\Theta,
 \qquad u\in U_v(\mathfrak g),
\end{equation}
where the bar on $U_v(\mathfrak g)\otimes U_v(\mathfrak g)$ acts on
each factor.  For based modules $X,Y$, the tensor-product bar is
\begin{equation}\label{eq:tensor-product-bar}
 \Psi_{X,Y}(x\otimes y)=\Theta(\bar x\otimes\bar y).
\end{equation}
These are the quasi-$R$-matrix and tensor bar of
\cite[Theorem~4.1.2 and Section~27.3]{Lusztig1993}, with tensor factors
interchanged to match \eqref{eq:coproduct-convention}
\cite{Lusztig1992Tensor}.  The sums act finitely on the finite-dimensional
modules used here.  In particular, a term with $\beta\ne0$ raises the
weight of the first factor by $\beta$ and lowers that of the second by
$\beta$.  This is the direction used in the minimum-vector bar checks.

An endpoint weight is \emph{squarefree} if every coordinate is zero
or one.  In degree $d$ on $d$ labels the squarefree weight is
\[
 \delta_d=\varepsilon_1+\cdots+\varepsilon_d.
\]
On this weight space write $\mathsf H_i=E_iF_i-v^{-1}\operatorname{id}$.
Sections~\ref{sec:row-model} and~\ref{sec:column-model} use the ascent
labels below in Gessel's expansion \eqref{eq:schur-fundamental} and
the bidirected connectivity of the slice for
condition~\ref{cond:levi-top-connectivity}.

\phantomsection\label{not:squarefree-tableaux}
Fix $\mu\vdash d$.  In the semistandard-tableau model of the type-$A$
crystal of $V(\mu)$ \cite{KashiwaraNakashima1994}, let $G_S$ denote
the lower global basis vector with crystal label $S$.  The weight of
$S$ records the number of occurrences of each letter.

For crystal calculations we read $S$ as the reverse of the row-reading
word $\operatorname{rr}(S)$, which reads each row from left to right,
starting with the bottom row (Subsection~\ref{subsec:insertion-molecules}).
For color $i$, retain only $i$ and $i+1$, write $+$ and $-$ for them,
respectively, and cancel adjacent $+-$ pairs until none remain.
The operator $\widetilde e_i$ changes the letter at the rightmost
uncancelled $-$ from $i+1$ to $i$, and $\widetilde f_i$ changes the
letter at the leftmost uncancelled $+$ from $i$ to $i+1$.
If the required sign is absent, the corresponding operator gives zero.
This is the signature rule for our coproduct and lower lattice at
$v=\infty$.

For a standard tableau $S$, recall from
Subsection~\ref{subsec:ascent-characteristics} that $\Des(S)$ consists of the $i$ for which $i+1$ lies in a lower row
than $i$.  On the squarefree slice we read the ascent label of $G_S$
from the Hecke action as in \eqref{eq:canonical-wgraph-action}.
It consists of the $i$ for which the coefficient of $G_S$ in
$\mathsf H_iG_S$ is $-v^{-1}$.

\begin{lemma}[Squarefree highest-weight slice]\label{lem:tableau-slice}
The crystal labels of $V(\mu)_{\delta_d}$ are the tableaux
$S\in\SYT(\mu)$.  For each such $S$, the coefficient of $G_S$ in
$\mathsf H_iG_S$ is
\[
 \begin{cases}
  -v^{-1},&i\in\Des(S),\\
  v,&i\notin\Des(S).
 \end{cases}
\]
In particular, the ascent label of $G_S$ is $\Des(S)$.
\end{lemma}
\begin{proof}
A tableau of shape $\mu$ has weight $\delta_d$ precisely when each of
$1,\ldots,d$ occurs once, so the tableaux of this weight are exactly
$\SYT(\mu)$.  If $i+1$ is below $i$, the reversed
row-reading word places $i$ before $i+1$, and the two signs cancel.
Otherwise they occur in the order $i+1,i$ and neither sign cancels.
Consequently
\[
 \varepsilon_i(S)=\varphi_i(S)=
 \begin{cases}0,&i\in\Des(S),\\1,&i\notin\Des(S).\end{cases}
\]
The rank-one crystal component therefore has highest weight $0$ or $2$,
respectively.  These are the only possibilities because the total content
in the two letters is two.

Apply Theorem~\ref{thm:excellent} to the rank-one Levi of color $i$.
In the associated graded module, the image of $G_S$ lies on the
zero-weight line of $V(0)$ or $V(2)$.  The normalized operator acts there by $-v^{-1}$ or $v$,
respectively.  The filtration and its quotients use subsets of the global
basis, so passage to the associated graded leaves the diagonal coefficient
unchanged.  Remark~\ref{rem:reduced-coefficient-example} gives an
example in which the generator also has an off-diagonal coefficient.
\end{proof}

\begin{corollary}[Bidirected connectivity of a lower squarefree layer]
\label{cor:lower-layer-connected}
For $\mu\vdash d$, the graph on the lower global basis of
$V(\mu)_{\delta_d}$ defined by the operators $\mathsf H_j$ is connected
after only bidirected edges are retained.
\end{corollary}
\begin{proof}
Fix $2\leq i\leq d-1$ and restrict the crystal to colors $i-1,i$.
Ignore all other letters in the reversed row-reading word and relabel
$i-1,i,i+1$ as $1,2,3$.  Applying the signature rule recalled before
Lemma~\ref{lem:tableau-slice} gives the following two paths, read from
left to right:
\[
 \begin{aligned}
 231&\xrightarrow{\widetilde e_1}131
     \xrightarrow{\widetilde e_2}121
     \xrightarrow{\widetilde f_1}122
     \xrightarrow{\widetilde f_2}132,\\
 312&\xrightarrow{\widetilde e_2}212
     \xrightarrow{\widetilde e_1}112
     \xrightarrow{\widetilde f_2}113
     \xrightarrow{\widetilde f_1}213.
 \end{aligned}
\]
The middle highest words $121$ and $112$ have weight $(2,1,0)$.
Normality identifies each of their components with $B(2,1)$, whose
squarefree weight has two vertices.  The tensor words $123$ and $321$
belong instead to $B(1,1,1)$ and $B(3)$, respectively, and each has a
one-point squarefree slice.  Reversing the words back to row-reading
order yields exactly
\begin{equation}\label{eq:rank-two-D-pairs}
 132\longleftrightarrow231,
 \qquad 213\longleftrightarrow312.
\end{equation}
By \eqref{eq:dual-equivalence}, these pairs are $S,D_iS$.
Letters outside the three-letter alphabet contribute no signs, so each
displayed path is also a path in the full tableau crystal.

By Theorem~\ref{thm:excellent} for this rank-two Levi, the associated
graded layer of top $(2,1)$ is a based direct sum of copies of
$V(2,1)$.  Consider the copy corresponding to the component just found.
Relabel its two simple roots as $1,2$.  As a
$U_v(\mathfrak{sl}_3)$-module it is $V(\omega_1+\omega_2)$.
Let $u$ be its highest global basis vector, normalized as the image of
the standard highest vector under the based identification, and put
\[
 e=F_2F_1u,\qquad f=F_1F_2u.
\]
These are the two global basis vectors of its squarefree weight.
Indeed, $F_1u=\widetilde f_1u$ is $2$-highest and
$F_2u=\widetilde f_2u$ is $1$-highest, because $E_2F_1u=E_1F_2u=0$.
It follows that
$e=\widetilde f_2\widetilde f_1u$ and
$f=\widetilde f_1\widetilde f_2u$ lie in the lower lattice and have
the two distinct squarefree crystal residues of $B(2,1)$: in the
standard tableau model, the highest tensor word is $112$, and these
two residues have tensor words $312$ and $213$, respectively.
Both vectors are integral and bar invariant, so balancedness identifies
them with the global lifts of those residues, without a scalar change.

Writing $[2]=v+v^{-1}$, the defining quantum-group relations give
\[
 E_1e=F_2u,\qquad E_1f=[2]F_2u,\qquad
 E_2e=[2]F_1u,\qquad E_2f=F_1u.
\]
The two vectors have zero coroot weight for both simple roots, so
$E_jF_j=F_jE_j$ on their span.  Hence
\[
 \mathsf H_1e=-v^{-1}e+f,\qquad \mathsf H_1f=vf,
 \qquad
 \mathsf H_2e=ve,\qquad \mathsf H_2f=e-v^{-1}f.
\]
In the original colors, $\mathsf H_{i-1}$ sends $e$ to a vector with
$f$-coefficient one, and $\mathsf H_i$ sends $f$ to a vector with
$e$-coefficient one.
The based filtration and its one-top decomposition preserve the
coefficients between the two surviving basis vectors.  Thus $S,D_iS$
are joined in both directions in the original squarefree graph.
Lemma~\ref{lem:dual-equivalence-connected} now proves connectivity.
\end{proof}

\section{\texorpdfstring{The row graph and its type-$C$ realization}
{The row graph and its type-C realization}}
\label{sec:row-model}
\label{subsec:row-positive-type-c}

We construct a lower based realization of each fixed-point block of
the row graph, verify
conditions~\ref{cond:levi-top-multiplicity}--\ref{cond:levi-top-connectivity},
and apply Theorem~\ref{thm:levi-top-criterion} to prove the row part of
Theorem~\ref{thm:main-intro}.  The construction assumes $d\geq2$;
the one-vertex case $d=1$ is included with the conclusions below.
We use the conventions of Subsection~\ref{subsec:squarefree-tools}.

\subsection{The based realization}

The radical generators form a quantum symmetric square, and the
fixed points contribute an exterior power.

\subsubsection{Root vectors and the based parabolic quotient}

Choose $R\geq\max\{2,d\}$.  For type $C_R$ take
\[
 \alpha_i=\varepsilon_i-\varepsilon_{i+1}\quad(1\leq i<R),
 \qquad \alpha_R=2\varepsilon_R,
 \qquad J_R=\{1,\ldots,R-1\}.
\]
The symmetrizing parameters are $v_i=v$ for $i<R$ and $v_R=v^2$.
The last fundamental weight is
$\varpi_R^C=\varepsilon_1+\cdots+\varepsilon_R$.
The longest element is $w_0^C=-1$, and $w_{0,J_R}$ denotes the
longest element of the Levi Weyl group.  Put
$w_C^J=w_0^Cw_{0,J_R}$ and, for $0\leq t<R$, put
\[
 c_t=s_Rs_{R-1}\cdots s_{t+1}.
\]
The expression
\begin{equation}\label{eq:c-reduced-word}
 w_C^J=c_0c_1\cdots c_{R-1}
       =s_{i_1}\cdots s_{i_m},\qquad m=R(R+1)/2,
\end{equation}
is reduced.  Its root sequence, defined by
\[
 \beta_s=s_{i_1}\cdots s_{i_{s-1}}(\alpha_{i_s})
 \qquad(1\leq s\leq m),
\]
lists
\begin{equation}\label{eq:c-radical-roots}
 \beta_{ab}=\varepsilon_a+\varepsilon_b
 \quad(1\leq a\leq b\leq R)
\end{equation}
by decreasing $b$ and then decreasing $a$; in particular,
$\beta_{aa}=2\varepsilon_a$.  To verify this, set
$P_t=c_0\cdots c_{t-1}$.  Induction gives
\[
 P_t(\varepsilon_a)=
 \begin{cases}
 -\varepsilon_{R+1-a},&a\leq t,\\
 \varepsilon_{a-t},&a>t.
 \end{cases}
\]
The first letter of block $c_t$ therefore contributes
$2\varepsilon_{R-t}$, and its subsequent letter $i$ contributes
$\varepsilon_{i-t}+\varepsilon_{R-t}$ for
$i=R-1,\ldots,t+1$.  These are precisely
$(b,b),(b-1,b),\ldots,(1,b)$ with $b=R-t$.
Finally $P_R(\varepsilon_a)=-\varepsilon_{R+1-a}$.
All roots in this sequence are positive and distinct, proving both
reducedness and the asserted Weyl-group identity.  For $R=2$ the word
is $(2,1,2)$ and its roots are $(22),(12),(11)$.

Recall the convention $\mathsf T_i=T'_{i,-1}$ from
\eqref{eq:shared-braid-convention}.  In this root order put
\[
 S_{ab}^{(r)}=
 \mathsf T_{i_1}\cdots\mathsf T_{i_{s-1}}(E_{i_s}^{(r)})
 \quad\text{when }\beta_s=\beta_{ab},\qquad S_{ab}=S_{ab}^{(1)}.
\]
Let $B_{C,R}^+$ be the quantum Schubert subalgebra generated by the
$S_{ab}$.  Its divided PBW products are the forward products in the
root order \eqref{eq:c-radical-roots}.  Put
\[
 \mathscr S_R=(B_{C,R}^+)^{\mathrm{op}},\qquad
 x\star y=yx\quad\text{in }B_{C,R}^+.
\]
Thus written $\star$-products have increasing right endpoints.
Each radical root has coefficient one at $\alpha_R$, and this
coefficient defines the degree.  The subspace supported on the first
$d$ endpoints is denoted by $\mathscr S_d$.

\begin{lemma}[The type-$C$ based parabolic quotient]
\label{lem:c-based-parabolic}
Let
\[
 Q^C_{R,\A}=
 U_\A^+(C_R)\Big/
 \operatorname{Sat}_{\A}
 \left(\sum_{j\in J_R}U_\A^+(C_R)E_j\right).
\]
The rational PBW map $B_{C,R}^+\longrightarrow
\mathbb K\otimes_\A Q^C_{R,\A}$ transports a lower based structure
to every finite homogeneous part of $\mathscr S_R$.
The divided radical PBW products are an $\A$-basis of its integral
lattice and an $\Oinf$-basis of its lower crystal lattice, and their
residues are its crystal basis.  A nontrivial term in the quotient
bar of a PBW product has strictly larger exponent data in the left
lexicographic order for \eqref{eq:c-radical-roots}.

For each fixed degree and all sufficiently large $N$, the map
\begin{equation}\label{eq:c-finite-embedding}
 [u]\longmapsto u v_{-N\varpi_R^C}
 \quad\text{into }V_{C_R}(N\varpi_R^C)
\end{equation}
is a based embedding for the derived type-$A$ Levi.  The ambient weight of $uv_{-N\varpi_R^C}$ is
$\wt(u)-N\mathbf 1_R$.  On the degree-$p$ part, the polynomial
weight lift of total degree $2p$ restores $\wt(u)$.  In particular,
every homogeneous part is an integrable lower based module for the Levi.
\end{lemma}
\begin{proof}
Recall from Subsection~\ref{subsec:quotient-conventions} that
$\varepsilon_j^*(b)$ counts how many times $\widetilde e_j^*$ can be
applied to $b$ before the result is zero.  By the ideal compatibility in the proof of
Lemma~\ref{lem:extremal-quotient}, the saturated ideal in the statement
is spanned by the $G(b)$ with $\varepsilon_j^*(b)>0$ for some $j\in J_R$.
Its complement therefore gives a based quotient.

Extend \eqref{eq:c-reduced-word} using
$w_0^C=w_C^Jw_{0,J_R}$.  Since
$w_C^J(\Phi_{J_R}^+)=\Phi_{J_R}^+$, the full PBW factorization has
the radical factors first and the Levi factors last.  Its nonzero
Levi-suffix data span the rational ideal
$\sum_{j\in J_R}U_v^+E_j$.  Intersecting this coordinate span with
the integral PBW lattice and with the lower PBW lattice retains
exactly the same suffix data.  The quotient consequently has the
asserted two PBW lattice bases and the surviving crystal residues.
Apply \eqref{eq:full-pbw-ordinary-bar} to the full word before
projecting.  Recall that the order in Lemma~\ref{lem:full-pbw-bar}
compares the first differing exponent from the left.  For surviving
tuples the Levi suffix is zero, so this comparison takes place in
the radical prefix.  All divided powers in this argument use their root
parameters, in particular $v_R=v^2$ at the long simple root.

Give the quotient its induced Levi action: $E_j$ acts by left
multiplication, $K_j$ by the weight, and $F_j$ by commuting through
a positive word to the Levi-trivial cyclic vector.  This action
preserves the $\alpha_R$-degree.  For $\Lambda=N\varpi_R^C$, the map $[u]\mapsto uv_{-N\varpi_R^C}$
is induced by $\pi_\Lambda$ of Subsection~\ref{subsec:quotient-conventions}.
Its bounds $m_i=\langle h_i,-w_0\Lambda\rangle$ are $m_j=0$ for $j\in J_R$
and $m_R=N$.  After the parabolic quotient, its kernel
is therefore the subset of surviving global labels with
$\varepsilon_R^*>N$.  There are finitely many PBW data in each
degree.  Taking $N$ at least the maximum of these finitely many
string parameters makes every label in that degree survive, with
coefficient one, by Lemma~\ref{lem:extremal-quotient}.
The lowest vector is trivial for the derived Levi, so the map is
linear for that algebra.  Its image is Levi-stable and is spanned
by a subset of a global basis.  Pullback proves integrability and the claimed lower
based structure, independently of sufficiently large $N$.
\end{proof}

\subsubsection{Endpoint action and straightening}

\begin{lemma}[Type-$C$ endpoint formulas]
\label{lem:c-endpoint-action}
On $B_{C,R}^+$ use the opposite-adjoint action
\[
 E_i\cdot a=E_i a-(K_i aK_i^{-1})E_i,\qquad
 F_i\cdot a=(F_i a-aF_i)K_i,\qquad
 K_i\cdot a=K_i aK_i^{-1}\quad(i\in J_R).
\]
It preserves the algebra and agrees with the induced action on
$Q^C_{R,\A}$.  Its degree-one module is
$V(2\omega_1)=\operatorname{Sym}_v^2(V_R)$, with lower global basis
$\{S_{ab}:a\leq b\}$.  Replacing an endpoint by an adjacent one
has coefficient one when all endpoints involved are distinct.
The remaining nonzero formulas, for $j=i+1$, are
\begin{equation}\label{eq:c-diagonal-string}
 \begin{aligned}
 E_iS_{jj}&=S_{ij},& E_iS_{ij}&=[2]S_{ii},\\
 F_iS_{ii}&=S_{ij},& F_iS_{ij}&=[2]S_{jj},
 \qquad [2]=v+v^{-1}.
 \end{aligned}
\end{equation}
The other endpoint actions are zero.  On $\mathscr S_R$ the product
rules are
\begin{equation}\label{eq:c-star-leibniz}
 \begin{aligned}
 E_i(a\star b)&=(E_i a)\star(K_i b)+a\star(E_i b),\\
 F_i(a\star b)&=(F_i a)\star b+(K_i^{-1}a)\star(F_i b).
 \end{aligned}
\end{equation}
\end{lemma}
\begin{proof}
Write $D_i$ for the displayed opposite-adjoint $E_i$-operator.
We first verify the nonsimply laced normalization.  In a rank-two
$C_2$ subsystem let $i$ be short and $j$ long, so
$a_{ij}=-2$, $a_{ji}=-1$, $v_i=v$, and $v_j=v^2$.
Put
\[
 r=v^{-2},\qquad
 D=E_iE_j-rE_jE_i=\mathsf T_j(E_i),\qquad
 L=-K_j^{-1}F_j=\mathsf T_j(E_j).
\]
The defining relations give
\[
 [F_j,D]=-rE_iK_j,\qquad LD=rDL+E_i.
\]
Applying $\mathsf T_j$ to
\[
 \mathsf T_i(E_j)
 =E_jE_i^{(2)}-v^{-1}E_iE_jE_i+v^{-2}E_i^{(2)}E_j
\]
and using $LD=rDL+E_i$ gives
\begin{align*}
 \mathsf T_j\mathsf T_i(E_j)
 &=\frac{LD^2}{[2]}-v^{-1}DLD+
                  r\frac{D^2L}{[2]}\\
 &=\frac{E_iD-DE_i}{[2]}\\
 &=E_i^{(2)}E_j-v^{-1}E_iE_jE_i+v^{-2}E_jE_i^{(2)}
   =D_i^{(2)}(E_j).
\end{align*}
Indeed, the coefficient of $D^2L$ is
$r(r+1)/[2]-v^{-1}r=0$.
The notation $D_i^{(2)}$ denotes the action of the divided power
$E_i^{(2)}$ in the opposite-adjoint representation.

These identities and the adjacent short-root braid relation give
the following recursions for the specified PBW vectors:
\begin{equation}\label{eq:c-root-recursion}
 S_{RR}=E_R,\qquad
 S_{ab}=D_aS_{a+1,b}\ (a<b),\qquad
 S_{aa}=D_a^{(2)}S_{a+1,a+1}\ (a<R).
\end{equation}
For a blockwise verification, set $P_t=c_0\cdots c_{t-1}$ and
$b=R-t$.  The adjacent short-root identities show inductively that
\[
 \mathsf T_{P_t}(E_j)=E_{j-t}
 \qquad(t+1\leq j\leq R-1).
\]
The first root in block $t$ is the diagonal root $(b,b)$.
The second is obtained by applying the first $C_2$ identity to
$\mathsf T_R(E_{R-1})$, and hence equals $D_{b-1}S_{bb}$.
The remaining roots in the block follow successively by the
adjacent short-root formula, giving $S_{ab}=D_aS_{a+1,b}$.
To pass to the first root of the next block, all letters of $c_t$
below $R-1$ commute with $E_R$.  The second $C_2$ identity above,
transported through $\mathsf T_{P_t}$, gives
$S_{b-1,b-1}=D_{b-1}^{(2)}S_{bb}$.  This proves
\eqref{eq:c-root-recursion} with coefficient one.

Under the opposite-adjoint action, $E_R$ is a lowest Levi-weight
vector of weight $2\varepsilon_R$.  All $F_i$ kill it, all
$E_i$ except $E_{R-1}$ kill it, and $E_{R-1}^3$ kills it by the
quantum Serre relation \eqref{eq:quantum-serre-relations}.
These are the lowest-weight relations of
$V(2\omega_1)$.  The recursions give its nonzero weight vectors,
one for each $\varepsilon_a+\varepsilon_b$ with $a\leq b$.
They thus identify the generated module with that simple module.
On its two-element strings the nonzero raising and lowering
coefficients are one.  On its three-element string, the
divided-power normalization in \eqref{eq:c-root-recursion} gives
exactly \eqref{eq:c-diagonal-string}.

The opposite-adjoint action is an action for
$\Delta^{\mathrm{op}}$ on ordinary multiplication; reversing the
product gives \eqref{eq:c-star-leibniz}.  It preserves the
algebra because it preserves the span of the generators.
On a Levi-trivial cyclic vector the rightmost $E_i,F_i$ terms
vanish and $K_i$ acts by one.  Therefore
$(u\cdot a)v_0=u(av_0)$, proving agreement with the quotient action.
Finally, in degree one every radical weight has a single surviving
PBW datum.  Quotient PBW bar triangularity fixes that vector, and
Lemma~\ref{lem:c-based-parabolic} identifies its crystal residue.
It is consequently the corresponding lower global vector.
\end{proof}

The endpoint formulas make the subspace supported on $[d]$ stable
under $J=\{1,\ldots,d-1\}$.  PBW straightening and nonnegative
endpoint weights make it a subalgebra.  Homogeneity gives the
restricted based structure on each $\mathscr S_d^{(p)}$.

\begin{lemma}[Two type-$C$ straightening relations]
\label{lem:c-two-relations}
In $\mathscr S_d$ one has
\begin{equation}\label{eq:c-two-relations}
 S_{bc}\star S_{ac}=vS_{ac}\star S_{bc}\quad(a<b<c),
 \qquad
 S_{ad}\star S_{bc}=S_{bc}\star S_{ad}\quad(a<b<c<d).
\end{equation}
\end{lemma}
\begin{proof}
Use the Levendorskii--Soibelman relation in the forward root order
and then reverse multiplication.  For the first pair the root inner
product is one.  A correction has weight
$\varepsilon_a+\varepsilon_b+2\varepsilon_c$ and must use roots
strictly between $(b,c)$ and $(a,c)$.
The two possible multisets of radical factors of that weight are
$\{(a,c),(b,c)\}$ and $\{(a,b),(c,c)\}$.
The first uses the endpoints of the interval; in the second,
$(c,c)$ lies before the interval and $(a,b)$ after it.
Neither is a correction supported strictly inside the interval.
The leading LS coefficient is therefore $v$.

For the second pair the inner product is zero.  Its four endpoint
coordinates are one, so no diagonal root can occur in a correction.
The alternative pairings are $\{(a,b),(c,d)\}$ and
$\{(a,c),(b,d)\}$.  In each case one root is before the interval
between $(a,d)$ and $(b,c)$ and the other is after it.
There is no correction, and the leading coefficient is one.
\end{proof}

\subsubsection{The squarefree tensor realization}

Let $\xi_F$ be the lower global vector in
$\bigwedge_v^k V_d$ corresponding to the $k$-element subset
$F\subseteq[d]$; for $k=0$ this denotes the trivial module.
When exactly one of $i,i+1$ lies in $F$, the nonzero $E_i$ or $F_i$
replaces it by the other with coefficient one.  When both or neither
lie in $F$, both operators kill $\xi_F$.
Put
\begin{equation}\label{eq:row-positive-tensor-model}
 \mathscr Y_{d;p,k}=\mathscr S_d^{(p)}\otimes
                         \bigwedge_v^k V_d,
 \qquad d=2p+k.
\end{equation}
Equip this tensor product with the bar of \eqref{eq:tensor-product-bar}.
Recall from Subsection~\ref{subsec:squarefree-tools} that at weight
$\delta_d$ each label occurs exactly once across the two factors.
The standard vectors of this weight are
\begin{equation}\label{eq:row-positive-standard-vectors}
 Z_w^{\mathrm r}
 =\left(\mathop{\star}_{(a,b)\text{ a two-cycle of }w}S_{ab}\right)
       \otimes\xi_{\Fix(w)},
 \qquad w\in I_d,\quad\#\Fix(w)=k,
\end{equation}
with the written pair factors in increasing right-endpoint order.
Squarefreeness forces every diagonal PBW exponent to be zero and
forces the pair endpoints and the singleton subset to be disjoint.
Thus these vectors are exactly the tensor PBW data of weight
$\delta_d$; they are an integral basis and a lower-lattice basis
of that slice, with the corresponding crystal residues.

\begin{theorem}[The type-$C$ row based comparison]
\label{thm:row-positive-comparison}
The squarefree Hecke operators of Subsection~\ref{subsec:squarefree-tools},
here denoted by
\begin{equation}\label{eq:row-positive-hecke}
 \mathsf H_i^{\mathrm r}=E_iF_i-v^{-1}\operatorname{id}
 \quad\text{on }(\mathscr Y_{d;p,k})_{\delta_d},
\end{equation}
give the row standard action \eqref{eq:four-cases-row} in
\eqref{eq:row-positive-standard-vectors}.  The assignment
$M_w\mapsto Z_w^{\mathrm r}$ identifies the integral lattices,
the standard bars, and the row canonical basis with the lower
global basis of this squarefree slice.
\end{theorem}
\begin{proof}
All endpoint weights are nonnegative.  On a squarefree weight the
two coordinates at $i,j=i+1$ have total two, so the local
zero-weight strings have highest weight zero or two.  The
operator in \eqref{eq:row-positive-hecke} has eigenvalues
$-v^{-1}$ and $v$ on those two lines.

First suppose $i,j$ belong to distinct pair carriers.  Write $A$
for the lower matching and $B$ for its conjugate by $s_i$.
Their affected factors and the weight-two target are
\[
\begin{array}{c|c|c|c}
 &A&B&Q\\\hline
 a<b<i&S_{ai}\star S_{bj}&S_{bi}\star S_{aj}&S_{ai}\star S_{bi}\\
 j<a<b&S_{ia}\star S_{jb}&S_{ja}\star S_{ib}&S_{ia}\star S_{ib}\\
 a<i<j<b&S_{ai}\star S_{jb}&S_{aj}\star S_{ib}&S_{ai}\star S_{ib}.
\end{array}
\]
All indices in these rows are distinct.  The coefficient-one
endpoint formulas, \eqref{eq:c-star-leibniz}, and
\eqref{eq:c-two-relations} give
\begin{equation}\label{eq:row-positive-carrier-ratio}
 E_iA=Q,\qquad E_iB=vQ,\qquad F_iQ=B+v^{-1}A.
\end{equation}
For instance, in the first row $E_iB=S_{bi}\star S_{ai}=vQ$,
and the first term of $F_iQ$ is
$S_{aj}\star S_{bi}=S_{bi}\star S_{aj}=B$.
In the other rows the factor $v$ in $E_iB$ comes from the
$K_i$-weight of the second affected factor.
Every unaffected factor is killed by $E_i,F_i$ and has
$K_i$-weight one, so it stays in its original place.  The only
changed right endpoints that can cross are $i,j$, and no third
right endpoint lies between them.  Consequently the displayed
straightening also applies inside a full matching.
Since $E_iF_i=F_iE_i$ on weight zero,
\eqref{eq:row-positive-carrier-ratio} yields
\[
 \mathsf H_i^{\mathrm r}A=B,\qquad
 \mathsf H_i^{\mathrm r}B=A+(v-v^{-1})B.
\]

For a pair--singleton strict ascent the affected tensors are
\[
 A=x_i\otimes\xi_{F\cup\{j\}},\qquad
 B=x_j\otimes\xi_{F\cup\{i\}},\qquad
 Q=x_i\otimes\xi_{F\cup\{i\}},
\]
where $F$ omits $i,j$ and the other endpoint of the pair is
unchanged.  The partner of $i$ in $A$ is an original label and
the completed partner of its singleton $j$ is an auxiliary label;
therefore this is the strict ascent in \eqref{eq:four-cases-row}.
The tensor coproduct and the exterior-power coefficients give
exactly \eqref{eq:row-positive-carrier-ratio}, hence the same
strict matrix.
If $(i,j)$ is a pair, \eqref{eq:c-diagonal-string} gives
$E_iF_iS_{ij}=[2]S_{ij}$ and hence eigenvalue $v$.
If both are singletons, $E_i,F_i$ kill the exterior carrier, and
the eigenvalue is $-v^{-1}$.  These are precisely the four row
standard cases.  In particular the Hecke relations hold, since
the resulting matrices are the prescribed standard action.

We next check the bar normalization.  Put $r=2p$.
Among the squarefree PBW data on $[r]$, the matching
$(12)(34)\cdots(r-1,r)$ is lexicographically maximal in the
forward root order \eqref{eq:c-radical-roots}.
All diagonal exponents are zero.  The first allowed pair is
$(r-1,r)$, and a matching containing it has exponent one there;
after deleting it, the same argument applies on $[r-2]$.
Thus Lemma~\ref{lem:c-based-parabolic} leaves no surviving
nontrivial bar correction to its PBW vector.
The minimum standard tensor is consequently the tensor of two
bar-fixed vectors,
\[
 (S_{12}\star S_{34}\star\cdots\star S_{r-1,r})
     \otimes\xi_{\{r+1,\ldots,d\}}.
\]
For $r=0$, the first factor is the empty product $1$.
For completeness, every nonconstant quasi-$R$ term kills this
tensor.  Its weight has the form
$\beta=\sum c_i\alpha_i\ne0$, with $c_i\geq0$, acting by a
positive-weight operator on the first factor and a
negative-weight operator on the second.
Let $s$ be the largest index with $c_s>0$.
If $s\geq r$, the first resulting weight has coordinate
$-c_s$ at $s+1$, outside its polynomial weight set.
If $s<r$, let $t$ be the smallest index with $c_t>0$; the second
resulting weight has coordinate $-c_t$ at $t$, which is also
impossible.  This includes the cases $p=0$ and $k=0$.
Thus \eqref{eq:tensor-product-bar} fixes the minimum.

Since the module bar fixes $E_i,F_i$,
\[
 \overline{\mathsf H_i^{\mathrm r}x}
 =(E_iF_i-v)\bar x
 =(\mathsf H_i^{\mathrm r})^{-1}\bar x.
\]
This bar and the row standard bar are compatible with the same
four-case matrices and fix the same minimum.  The strict-ascent
argument in Subsection~\ref{subsec:reduced-graphs} therefore identifies
the two bars.  The standard integral lattices and their lower crystal
lattices, including their residues, were identified in
\eqref{eq:row-positive-standard-vectors}.
By \eqref{eq:canonical-normalization}, the row canonical vectors are
bar fixed and congruent to the standard vectors modulo $v^{-1}$
times this lower lattice.  Balanced uniqueness
\eqref{eq:balanced-uniqueness-prelim} identifies them with the lower
global vectors, proving the based comparison.
\end{proof}

\subsection{Verification of the three conditions}

\begin{lemma}[Multiplicity-free character of the row realization]
\label{lem:row-positive-character}
The type-$A_{d-1}$ character of \eqref{eq:row-positive-tensor-model}
is
\begin{equation}\label{eq:row-positive-pieri}
 \operatorname{ch}\mathscr Y_{d;p,k}
 =\sum_{\substack{\mu\vdash d\\
             \mu\text{ has exactly }k\text{ odd rows}}}s_\mu(X).
\end{equation}
Every summand occurs with multiplicity one.
\end{lemma}
\begin{proof}
The independent exponents of the radical PBW generators give
\[
 \sum_{p\geq0}t^p\operatorname{ch}\mathscr S_d^{(p)}
 =\prod_{1\leq a\leq b\leq d}(1-tX_aX_b)^{-1}
 =\sum_{\lambda:\,\lambda\text{ has even parts}}
       t^{|\lambda|/2}s_\lambda(X).
\]
The final equality is Littlewood's symmetric-square identity
\cite[Chapter~I, \S5]{Macdonald1995}.  It applies to formal power
series, so extracting each degree is legitimate.
The exterior factor has character $e_k=s_{(1^k)}$.
By the vertical-strip Pieri rule, the coefficient of $s_\mu$
in its product with the degree-$p$ term counts even-row
partitions $\lambda$ for which $\mu/\lambda$ is a vertical
$k$-strip.
If such a $\lambda$ exists, a box is removed from exactly
each odd row of $\mu$, and from no even row.  This determines
$\lambda$ uniquely.  Conversely, replacing each part $\mu_i$
by $2\lfloor\mu_i/2\rfloor$ gives a partition with even parts;
the removed boxes form a vertical strip whose size is the
number of odd rows of $\mu$.
It has degree $2p$ exactly when that number is $k$.
This proves \eqref{eq:row-positive-pieri} and multiplicity one.
\end{proof}

For a row vertex $w$ in the fixed-$k$ block, let
$\operatorname{top}^{C}_J(w)$ be its lower Levi top in
$\mathscr Y_{d;p,k}$ under Theorem~\ref{thm:row-positive-comparison}.
Put
\[
 \Lambda^{\mathrm r}_{d;k}
 =\{\mu\vdash d:\mu\text{ has exactly }k\text{ odd rows}\}.
\]

\begin{proposition}[The three conditions for the row model]
\label{prop:row-criterion-hypotheses}
Fix $d=2p+k\geq2$.  The fixed-$k$ row block satisfies conditions
\ref{cond:levi-top-multiplicity}--\ref{cond:levi-top-connectivity}
for $\tau(w)=\operatorname{top}^{C}_J(w)$ and
$\Lambda=\Lambda^{\mathrm r}_{d;k}$, ordered by dominance.
For every dominance upset $\Omega\subseteq\Lambda^{\mathrm r}_{d;k}$,
the space
\begin{equation}\label{eq:row-positive-top-upset}
 \mathscr M^{\mathrm r}_\Omega
 =\bigoplus_{\operatorname{top}^{C}_J(w)\in\Omega}
       \A\,\underline M_w
\end{equation}
is a saturated $\Hh_d$-submodule.  For each occurring top $\mu$,
the one-top subquotient is based-isomorphic to the squarefree lower
slice $V(\mu)_{\delta_d}$, and its top fiber is bidirected connected
in the original row graph.
\end{proposition}
\begin{proof}
For condition~\ref{cond:levi-top-multiplicity},
Lemma~\ref{lem:row-positive-character} shows that exactly the weights
in $\Lambda^{\mathrm r}_{d;k}$ occur and that each occurs once.
Recall that the \hyperref[not:squarefree-tableaux]{tableau label $S$ of
$G_S$} has weight $\delta_d$ when each letter occurs once.
Lemma~\ref{lem:tableau-slice} thus indexes a basis of
$V(\mu)_{\delta_d}$ by the nonempty set $\SYT(\mu)$.
The based comparison in Theorem~\ref{thm:row-positive-comparison}
therefore gives the required surjective multiplicity-free top datum.

For condition~\ref{cond:levi-top-upsets}, apply
Theorem~\ref{thm:excellent} to the full finite based module
$\mathscr Y_{d;p,k}$.  It gives a saturated Levi submodule spanned
by the global vectors whose tops lie in $\Omega$.  The operator
$\mathsf H_i^{\mathrm r}=E_iF_i-v^{-1}\operatorname{id}$ preserves
this submodule and has weight zero, so it preserves its squarefree
slice.  Theorem~\ref{thm:row-positive-comparison} transports that
slice to exactly \eqref{eq:row-positive-top-upset}.  The complementary
canonical vectors form a free basis of the quotient, proving saturation.

For condition~\ref{cond:levi-top-connectivity}, fix
$\mu\in\Lambda^{\mathrm r}_{d;k}$ and take the upsets
\[
 \Omega_1=\{\nu\in\Lambda^{\mathrm r}_{d;k}:\mu\preceq\nu\},
 \qquad
 \Omega_0=\Omega_1\setminus\{\mu\}.
\]
Theorem~\ref{thm:excellent} and multiplicity one identify
$\mathscr M^{\mathrm r}_{\Omega_1}/
 \mathscr M^{\mathrm r}_{\Omega_0}$ with the squarefree lower based
slice of the unique copy of $V(\mu)$.
Its bidirected graph is connected by
Corollary~\ref{cor:lower-layer-connected}.
The filtration terms and the quotient are based on actual subsets
of the original canonical basis.  Taking the quotient deletes only
the vectors with top strictly above $\mu$ and leaves unchanged every
coefficient between vectors of top $\mu$.  Hence each bidirected edge
in the slice remains a bidirected edge in the original row graph.
This proves connectivity of the whole top fiber.
\end{proof}

\subsection{Cells, insertion shapes, and dominance}
\label{subsec:row-results}

Proposition~\ref{prop:row-criterion-hypotheses} verifies
conditions~\ref{cond:levi-top-multiplicity}--\ref{cond:levi-top-connectivity}
for the row realization constructed in this section.  We now apply
Theorem~\ref{thm:levi-top-criterion} and identify its top fibers by
row insertion.  The conclusions below include $d=1$: its graph has
one vertex, with polynomial top $(1)$, and all assertions are
immediate in that case.  For $d\geq2$, fix the number $k$ of fixed
points and put $p=(d-k)/2$.  The lower Levi top in the type-$C$
model is $\operatorname{top}^{C}_J$.

\begin{theorem}[$\mcell$-cells, $\mcell$-molecules, and row Levi tops]
\label{thm:m-cells-molecules}
In each block with $k$ fixed points, the $\mcell$-cells and
$\mcell$-molecules are precisely the lower Levi-top fibers
\[
 C^{\mathrm r}_\mu
 =\{w\in I_d:\#\Fix(w)=k,\quad
             \operatorname{top}^{C}_J(w)=\mu\},
\]
where $\mu$ ranges over the occurring tops.  In particular, every
$\mcell$-cell of $\mGamma$ is exactly one $\mcell$-molecule.
\end{theorem}
\begin{proof}
For $d\geq2$, Proposition~\ref{prop:row-criterion-hypotheses}
verifies all three hypotheses of Theorem~\ref{thm:levi-top-criterion}
in each fixed-$k$ block.  The criterion therefore identifies the
cells, molecules, and the
fibers $C^{\mathrm r}_\mu$.  By
Subsection~\ref{subsec:reduced-graphs}, no coefficient arrow crosses
fixed-point blocks, so the assertion also holds in the full row graph.
\end{proof}

\subsubsection{The row Levi label}

For an occurring top $\mu$ in a fixed-$k$ block, let
$C^{\mathrm r}_\mu$ be the fiber defined in
Theorem~\ref{thm:m-cells-molecules}.
By Proposition~\ref{prop:row-criterion-hypotheses}, its one-top
associated graded is the squarefree lower based slice of $V(\mu)$.
Recall the \hyperref[not:squarefree-tableaux]{diagonal-coefficient
convention} before Lemma~\ref{lem:tableau-slice}.  If the image of
$\underline M_w$ in this layer corresponds to $G_S$, the diagonal
coefficient of $H_i$ is preserved in the quotient.  The lemma therefore
gives $\mAsc(w)=\Des(S)$.  Summing over the layer and applying
\eqref{eq:schur-fundamental} yields
\begin{equation}\label{eq:row-top-character}
 \sum_{w\in C^{\mathrm r}_\mu}F_{\mAsc(w)}
 =\sum_{S\in\SYT(\mu)}F_{\Des(S)}=s_\mu.
\end{equation}
By Theorem~\ref{thm:m-cells-molecules}, this fiber is one row
molecule.  Proposition~\ref{prop:reduced-molecules} identifies it
with an insertion-shape fiber of some shape $\lambda$.
Comparing \eqref{eq:row-top-character} with the row identity
\eqref{eq:row-shape-character-new} in
Proposition~\ref{prop:insertion-shape-characteristics} yields
$s_\mu=s_{\lambda^{\tr}}$.  Linear independence of Schur functions
therefore proves
\begin{equation}\label{eq:row-label-identification}
 \operatorname{top}^{C}_J(w)=\lambda_{\mathrm r}(w)^{\tr}.
\end{equation}

\subsubsection{Row dominance and the cells}

\begin{theorem}[Row dominance filtration]
\label{thm:row-dominance-filtration}
For every dominance downset $\mathscr D$ of partitions of $d$,
\begin{equation}\label{eq:row-downset-module}
 \mathcal F^{\mathrm r}_{\mathscr D}
 =\bigoplus_{\lambda_{\mathrm r}(w)\in\mathscr D}
       \A\,\underline M_w
\end{equation}
is a saturated $\Hh_d$-submodule of the row model.
\end{theorem}
\begin{proof}
In each fixed-$k$ block, \eqref{eq:row-label-identification}
identifies the indicated top set with
$\{\mu:\mu^{\tr}\in\mathscr D\}$ among its occurring tops.
By \eqref{eq:transpose-reverses}, this is an upset.
Equation~\eqref{eq:row-positive-top-upset} makes its canonical-basis
span stable and saturated.  Taking the direct sum over all
admissible $k$ proves the assertion.  Saturation also follows
directly from the free complementary canonical-basis span.
\end{proof}

\begin{theorem}[Row edge dominance]
\label{thm:row-edge-dominance}
Let $y,z\in I_d$.  If
$\iota_{\mathrm{asc}}(y)\longrightarrow_{\mcell}
\iota_{\mathrm{asc}}(z)$ is a coefficient arrow in the reduced row
graph, then
\[
 \lambda_{\mathrm r}(z)\preceq\lambda_{\mathrm r}(y).
\]
If the reverse arrow also occurs, the two shapes are equal.
\end{theorem}
\begin{proof}
Apply Theorem~\ref{thm:row-dominance-filtration} to the principal
downset of $\lambda_{\mathrm r}(y)$.  It contains
$\underline M_y$, and its stability forces every nonzero
canonical-basis coefficient of $H_i\underline M_y$ to remain
in that subset span.  In particular,
$\lambda_{\mathrm r}(z)\preceq\lambda_{\mathrm r}(y)$.
Apply the same implication in both directions for a bidirected edge.
\end{proof}

\begin{corollary}[Row cells by insertion shape]
\label{cor:row-cells-shape}
Two vertices lie in the same $\mcell$-cell if and only if their
modified row Beissinger tableaux have the same shape.
\end{corollary}
\begin{proof}
Combine Theorem~\ref{thm:m-cells-molecules} with
Proposition~\ref{prop:reduced-molecules}.
\end{proof}

Corollary~\ref{cor:row-cells-shape} and
Theorem~\ref{thm:row-edge-dominance} prove
Theorem~\ref{thm:main-intro}\textup{(a)}.
Theorem~\ref{thm:row-dominance-filtration} proves the row filtration
asserted in the final sentence of Theorem~\ref{thm:main-intro}.

\newpage
\subsection{A rank-six row cell}

\begin{example}[A complete $\mcell$-cell $W$-graph in rank six]
\label{ex:rank-six-m-cell}
Work in the fixed-point-free block of $I_6$.  Let
\[
\begin{array}{c|c|c|c}
 x&w_x&\PrB(w_x)&\mAsc(w_x)\\ \hline
 A&(12)(34)(56)&\rowtableau{1&3&5\cr2&4&6\cr}&\{2,4\}\\
 B&(12)(35)(46)&\rowtableau{1&3&4\cr2&5&6\cr}&\{2,3,5\}\\
 C&(13)(24)(56)&\rowtableau{1&2&5\cr3&4&6\cr}&\{1,3,4\}\\
 D&(13)(25)(46)&\rowtableau{1&2&4\cr3&5&6\cr}&\{1,3,5\}\\
 E&(14)(25)(36)&\rowtableau{1&2&3\cr4&5&6\cr}&\{1,2,4,5\}.
\end{array}
\]
All five restricted row tableaux have shape $(3,3)$.  Therefore
Proposition~\ref{prop:reduced-molecules} and
Theorem~\ref{thm:m-cells-molecules} show that
\[
 \mathcal C_{\mcell}^{(3,3)}
 =\{w_A,w_B,w_C,w_D,w_E\}
\]
is one $\mcell$-cell and one $\mcell$-molecule.

Applying the standard matrices \eqref{eq:four-cases-row}, the bar
recursion from the block minimum, and the normalization
\eqref{eq:canonical-normalization} gives
the following complete reduced $W$-graph of its cell subquotient.  In the
diagram a solid double arrow is a bidirected edge and the dashed arrow is
one-directional.  The table following the diagram records every directed
coefficient: an entry $x\to y:J$ means that, for each $i\in J$, the target
canonical basis vector at $y$ occurs with coefficient $1$ in $H_i$ applied
to the source at $x$.  All off-diagonal weights are $1$, and no other
off-diagonal arrows occur.  Together with the ascent sets above, these data
record the entire reduced $W$-graph and its Hecke action.
\[
\begin{tikzpicture}[
  x=1.25cm,y=1.05cm,
  vertex/.style={circle,draw,thick,minimum size=7mm,inner sep=0pt},
  bi/.style={{Stealth[length=2.1mm]}-{Stealth[length=2.1mm]},semithick},
  one/.style={-{Stealth[length=2.1mm]},thick,dashed}
]
 \node[vertex] (A) at (0,2.5) {$A$};
 \node[vertex] (B) at (-2,1) {$B$};
 \node[vertex] (C) at (2,1) {$C$};
 \node[vertex] (D) at (0,-0.5) {$D$};
 \node[vertex] (E) at (0,-2.2) {$E$};

 \draw[bi] (A) -- (B);
 \draw[bi] (A) -- (C);
 \draw[bi] (B) -- (D);
 \draw[bi] (C) -- (D);
 \draw[bi] (D) -- (E);
 \draw[one,bend right=58] (E) to (A);
\end{tikzpicture}
\]
\[
\begin{array}{c|c@{\qquad}c|c@{\qquad}c|c}
\text{arrow}&J&\text{arrow}&J&\text{arrow}&J\\ \hline
A\to B&\{4\}&B\to A&\{3,5\}&A\to C&\{2\}\\
C\to A&\{1,3\}&B\to D&\{2\}&D\to B&\{1\}\\
C\to D&\{4\}&D\to C&\{5\}&D\to E&\{3\}\\
E\to D&\{2,4\}&E\to A&\{1,5\}&\multicolumn{2}{c}{}
\end{array}
\]
In particular, the five pairs joined in both directions form the complete
$\mcell$-molecule, while the long arrow $E\to A$ is an additional
one-directional edge inside the same cell.
\end{example}

\section{\texorpdfstring{The column graph and its type-$D$ realization}
{The column graph and its type-D realization}}
\label{sec:column-model}
\label{subsec:column-positive-type-d}

We construct the type-$D$ lower based realization of the column
blocks and verify
conditions~\ref{cond:levi-top-multiplicity}--\ref{cond:levi-top-connectivity}.
Applying Theorem~\ref{thm:levi-top-criterion} then proves the column
part of Theorem~\ref{thm:main-intro}.  The construction assumes
$d\geq2$, and the conclusions also include the one-vertex case
$d=1$.  We use the common tools and conventions of
Section~\ref{sec:levi-top-criterion}.

\subsection{The based realization}

The radical generators form a quantum exterior square, and the
fixed points contribute a symmetric power.

\subsubsection{Root vectors and the multiplication convention}

Choose an even integer
$R\geq\max\{4,d\}$.  (The harmless enlargement from $d$ to $R$ is useful
only for avoiding the parity interchange of the two half-spin weights.)  In
type $D_R$ choose simple roots
\[
 \alpha_i=\varepsilon_i-\varepsilon_{i+1}\quad(1\leq i<R),
 \qquad
 \alpha_R=\varepsilon_{R-1}+\varepsilon_R,
\]
where the first $R-1$ roots form a Levi subsystem $J_R$ of type
$A_{R-1}$.  The symmetrizing parameters are $v_i=v$ for all $i$.
Let $w_0$ and $w_{0,J_R}$ be the longest elements of the full
and Levi Weyl groups, respectively, and put $w^{J_R}=w_0w_{0,J_R}$; write
$\operatorname{Inv}(w)$ for the inversion set.  A direct root calculation gives
\begin{equation}\label{eq:radical-roots}
 \operatorname{Inv}(w^{J_R})=
 \{\varepsilon_a+\varepsilon_b:1\leq a<b\leq R\}.
\end{equation}
Every root in \eqref{eq:radical-roots} has coefficient one at the deleted spin
node.  Put \(m=R(R-1)/2\).  For \(0\leq t\leq R-2\), define
\[
 c_t=s_{R-(t\bmod 2)}s_{R-2}s_{R-3}\cdots s_{t+1},
\]
where the descending string following the first factor is omitted when
\(t+1>R-2\).  We use the reduced expression
\[
 w^{J_R}=c_0c_1\cdots c_{R-2}=s_{i_1}\cdots s_{i_m}.
\]
Its root sequence \(\beta_k=s_{i_1}\cdots s_{i_{k-1}}(\alpha_{i_k})\)
lists the roots \(\varepsilon_a+\varepsilon_b\) first by decreasing
\(b\), then by decreasing \(a\).  Indeed, induction on the blocks
\(c_t\) gives this root sequence; all its roots are positive, and the
resulting signed permutation sends \(\varepsilon_a\) to
\(-\varepsilon_{R+1-a}\).  Thus it is \(w_0w_{0,J_R}\), and the word
is reduced.  For example, when \(R=4\) the word is
\((4,2,1,3,2,4)\), with root sequence
\((34),(24),(14),(23),(13),(12)\).

With $\mathsf T_i=T'_{i,-1}$ as in \eqref{eq:shared-braid-convention},
the divided root vectors in this order are
\[
 Z_{\beta_a}^{(c)}=
 \mathsf T_{i_1}\cdots\mathsf T_{i_{a-1}}(E_{i_a}^{(c)}),
 \qquad Z_{\beta_a}=Z_{\beta_a}^{(1)},
 \qquad Z_{ut}=Z_{\varepsilon_u+\varepsilon_t}.
\]
Let \(B_R^+\) be the positive quantum Schubert subalgebra generated by
these vectors.  Its PBW products are the forward products
\begin{equation}\label{eq:correct-forward-pbw}
 P(\mathbf c)=Z_{\beta_1}^{(c_1)}\cdots Z_{\beta_m}^{(c_m)}.
\end{equation}
For matching notation we use the opposite algebra
\[
 \mathscr A_R=(B_R^+)^{\mathrm{op}},\qquad a\star b=ba
 \quad\text{in the underlying algebra }B_R^+.
\]
Thus \eqref{eq:correct-forward-pbw} is written in \(\mathscr A_R\) as
\(Z_{\beta_m}^{(c_m)}\star\cdots\star Z_{\beta_1}^{(c_1)}\).
The written \(\star\)-factors have increasing right endpoints, whereas
the ordinary product in \(B_R^+\) has decreasing right endpoints.
We use divided powers in all degrees.  In squarefree weight every exponent
is zero or one, so the divided and ordinary powers agree.

The algebra $\mathscr A_d$ is the subspace spanned by the PBW
products using endpoints in $\{1,\ldots,d\}$.  We use the initial
Levi $J$ and squarefree weight $\delta_d$ of
Subsection~\ref{subsec:squarefree-tools}.

\subsubsection{The based parabolic quotient}

Let $\varpi_R=\tfrac12\mathbf 1_R$ be the half-spin fundamental
weight corresponding to the deleted node $\alpha_R$.  Since $R$ is
even, the extremal weight $w^{J_R}(N\varpi_R)$ is $-N\varpi_R$.
Its extremal vector is trivial for the derived Levi
$U_v(\mathfrak g_{J_R})$.  Multiplication on that vector defines
\[
 \iota_N:B_R^+\longrightarrow V(N\varpi_R),
 \qquad x\longmapsto x v_{-N\varpi_R}.
\]

\begin{lemma}[Parabolic based realization]\label{lem:parabolic-based}
Let
\[
 I_{J_R}=
 \operatorname{Sat}_{\A}\!\left(\sum_{j\in J_R}U_{\A}^+E_j\right).
\]
The quotient \(Q_R=U_\A^+/I_{J_R}\), extended to \(\mathbb Q(v)\), is a
lower based graded vector space.  The length-additive factorization
\[
 w_0=w^{J_R}w_{0,J_R}
\]
and its PBW factorization identify \(Q_R\), as a graded vector space, with
\(B_R^+\), the underlying vector space of \(\mathscr A_R\).  This transports to every homogeneous part of
\(\mathscr A_R\) a lower global basis and a semilinear degreewise bar.
\end{lemma}
\begin{proof}
Recall from Subsection~\ref{subsec:quotient-conventions} that
$\varepsilon_j^*(b)>0$ precisely when $\widetilde e_j^*b\ne0$.
By the ideal compatibility
in the proof of Lemma~\ref{lem:extremal-quotient}, the saturation of
$U_\A^+E_j$ is spanned by the $G(b)$ with $\varepsilon_j^*(b)>0$.
Summing over $j\in J_R$ gives the span of the union of these subsets.  It is a based left ideal, so the complementary
global basis classes make \(Q_R\) a based quotient.

Since \(R\) is even,
\[
 w_0=-1,\qquad w^{J_R}=-w_{0,J_R},\qquad
 w^{J_R}(\Phi_{J_R}^+)=\Phi_{J_R}^+.
\]
A full PBW basis extending the radical root vectors gives the rational
vector-space factorization
\[
 U_v^+=B_R^+\,U_v^+(\mathfrak g_{J_R}).
\]
Equivalently, straightening moves the Levi factors to the right, with
nonzero powers of $v$ on the diagonal.  Nonconstant rightmost Levi
factors vanish in the quotient, and the radical PBW products give a
$\mathbb K$-basis of its surviving classes.  Hence
$x\mapsto x+I_{J_R}$ is a rational graded isomorphism.
The integral and crystal-lattice comparison is proved in
Lemma~\ref{lem:radical-pbw-balanced}.

We transport the quotient global basis, lower lattice and its bar image,
and bar through this identification.  The bar here is a module bar
transported from the quotient.  Homogeneity and nonnegativity of the
\(\varepsilon\)-weights imply that the basis elements supported on the first
\(d\) indices form the transported lower global basis of \(\mathscr A_d\).
\end{proof}

\begin{lemma}[Spin approximation of the based parabolic quotient]
\label{lem:spin-approx}
Give \(Q_R\) its grading by the coefficient of \(\alpha_R\).
Every homogeneous part is a finite lower based integrable
\(U_v(\mathfrak g_{J_R})\)-module.  The map
\[
 \bar\iota_N:Q_R\longrightarrow V(N\varpi_R),\qquad
 [u]\longmapsto u v_{-N\varpi_R},
\]
is, in every fixed degree and for all sufficiently large \(N\), a
based embedding for the derived Levi.  Its action on the source
is the action induced on the parabolic quotient.
The ambient weight of the image of a homogeneous $u$ is
$\wt(u)-\tfrac N2\mathbf 1_R$.  This constant shift vanishes on
the derived Levi coroots.  The polynomial lift of total degree
$2p$ assigns $\wt(u)$ to the image of a degree-$p$ vector.
\end{lemma}
\begin{proof}
The quotient is the positive-half realization of the parabolically
induced cyclic module whose cyclic vector is trivial for the Levi.
In particular, \(E_j\) acts by left multiplication, \(K_j\) by
the weight, and \(F_j\) acts by commuting it through a positive
word and then using \(F_j1=0\), \(K_j1=1\), for \(j\in J_R\).
The quantum-group relations make these actions well defined.
They preserve the \(\alpha_R\)-degree and do not depend on \(N\).

For $\Lambda=N\varpi_R$, the map $\iota_N$ is the restriction of
$\pi_\Lambda$ from Subsection~\ref{subsec:quotient-conventions}.
Recall that $m_i=\langle h_i,-w_0\Lambda\rangle$.  Here $w_0=-1$,
so $m_j=0$ for $j\in J_R$ and $m_R=N$.
The first conditions produce exactly the based quotient \(Q_R\).
The remaining kernel is spanned by those surviving global-basis
labels whose \(\varepsilon_R^*\) exceeds \(N\).
The PBW theorem makes each fixed-degree part finite dimensional.
Choose \(N\) at least the maximum of \(\varepsilon_R^*\) over
the finite set of labels in that degree.  The quotient map then
sends them to distinct global-basis vectors with coefficient one.
It is linear for the derived Levi because the cyclic vector of
the spin module is trivial for that algebra and hence satisfies
the relations defining the induced action.  Its image is
Levi-stable and spanned by a subset of a global basis.  Pulling back
the based structure proves the assertion and its independence of
sufficiently large \(N\).
\end{proof}

\begin{lemma}[The PBW lattice in the parabolic quotient]
\label{lem:radical-pbw-balanced}
Under the rational PBW identification
\(B_R^+\longrightarrow Q_R\), \(a\mapsto[a]\), the products
\eqref{eq:correct-forward-pbw} form an \(\A\)-basis of the integral
quotient and an \(\Oinf\)-basis of its lower crystal lattice.  Their
residues are the surviving crystal basis of the quotient.
In the quotient bar of a PBW product, each nontrivial term has
exponent data strictly larger in the ordinary left lexicographic order
on $(c_1,\ldots,c_m)$, with the total root weight fixed.
\end{lemma}
\begin{proof}
Extend the chosen reduced expression to
\(w_0=w^{J_R}w_{0,J_R}\), and use the same braid convention on the
resulting full PBW basis.  The radical factors precede the Levi factors.

For this full word, Lemma~\ref{lem:full-pbw-bar} and the lattice
calculation following its proof give the ordinary bar expansion,
the two PBW lattice bases and their lower crystal residues.

The full PBW factorization shows that
\(\sum_{j\in J_R}U_v^+E_j\) is the rational span of exactly those PBW
vectors having a nonzero Levi suffix.  Intersecting this coordinate span
with either of the two PBW lattices retains exactly the same subset of
basis vectors.  By the right-ideal compatibility in
Lemma~\ref{lem:extremal-quotient}, this is also a based ideal.
Passing to the quotient deletes the nonzero Levi-suffix data and leaves
precisely \eqref{eq:correct-forward-pbw}.  This proves all three basis
assertions.  In \eqref{eq:full-pbw-ordinary-bar}, any two surviving
data have zero Levi suffix.  In the left lexicographic order fixed before
Lemma~\ref{lem:full-pbw-bar}, their first differing exponent therefore
lies in the radical prefix.
\end{proof}

For homogeneous \(a\in B_R^+\), define the Levi action by the following
formulas, with products taken in the ambient type-\(D_R\) quantum group.
\begin{equation}\label{eq:opposite-adjoint-positive}
 \begin{aligned}
 E_i\cdot a&=E_i a-(K_i aK_i^{-1})E_i,\\
 F_i\cdot a&=(F_i a-aF_i)K_i,\\
 K_i\cdot a&=K_i aK_i^{-1}.
 \end{aligned}
\end{equation}
These are the opposite Hopf-adjoint formulas
\(u\cdot a=\sum_{(u)}u_{(2)}aS^{-1}(u_{(1)})\).
They are distinct from the ordinary adjoint action.

\subsubsection{Endpoint action and straightening}

\begin{lemma}[Endpoint formulas and the Levi action]
\label{lem:positive-endpoint-action}
The action \eqref{eq:opposite-adjoint-positive} preserves \(B_R^+\).
On its radical generators, \(E_i\) replaces an endpoint \(i+1\) by
\(i\), and \(F_i\) replaces an endpoint \(i\) by \(i+1\), with
coefficient one whenever the resulting endpoints are distinct.  Every
other action is zero.  In particular,
\[
 E_iZ_{a,i+1}=Z_{a,i},\qquad E_iZ_{i+1,b}=Z_{i,b},\qquad
 E_iZ_{i,i+1}=F_iZ_{i,i+1}=0.
\]
The action on \(B_R^+\) is a module-algebra action for
\(\Delta^{\mathrm{op}}\); on \(\mathscr A_R\) it is a
module-algebra action for \(\Delta\).  It agrees with the induced Levi
action on \(Q_R\).
\end{lemma}
\begin{proof}
Write \(D_i(a)=E_i\cdot a\).  An explicit recursion for the root
vectors is
\begin{equation}\label{eq:endpoint-root-recursion}
 Z_{aR}=D_aD_{a+1}\cdots D_{R-2}(E_R),\qquad
 Z_{ab}=D_bD_{b+1}\cdots D_{R-1}(Z_{aR})\quad(b<R).
\end{equation}
An empty string of operators is omitted.  Each nonzero step replaces one
endpoint; its current coroot weight is \(-1\), so the step is explicitly
\(a\mapsto E_i a-v^{-1}aE_i\).
Here is a blockwise check of \eqref{eq:endpoint-root-recursion} inside
\(U_v^+\), before passing to the quotient.  Put
\(\eta_t=R-(t\bmod2)\), \(P_t=c_0\cdots c_{t-1}\), and
\(b=R-t\).  The induction invariants are
\[
 \mathsf T_{P_t}(E_j)=E_{j-t}\quad(t+1\leq j\leq R-2),\qquad
 \mathsf T_{P_t}(E_{\eta_{t+1}})=E_{R-t-1}.
\]
They hold for \(t=0\).  For \(0\leq t\leq R-3\), the adjacent-node
identity \(\mathsf T_i\mathsf T_j(E_i)=E_j\), together with
commutation at nonadjacent nodes, gives
\[
 \mathsf T_{c_t}(E_j)=E_{j-1}\quad(t+2\leq j\leq R-2),\qquad
 \mathsf T_{c_t}(E_{\eta_t})=E_{R-2}.
\]
Composing with \(\mathsf T_{P_t}\) proves the two invariants for
\(t+1\).  In block \(t\), the first root is \((b-1,b)\).
For each subsequent root \((a,b)\), apply the adjacent braid formula
to the next letter \(a+t\); the preceding letters in that block
commute with \(E_{a+t}\), except for its immediate adjacent
predecessor.  The first invariant therefore gives
\[
 Z_{ab}=E_aZ_{a+1,b}-v^{-1}Z_{a+1,b}E_a.
\]
For the first root of a new block \(t\geq1\), all letters of
\(c_{t-1}\) below \(R-2\) commute with \(E_{\eta_t}\), and
the two spin nodes commute.  Applying the adjacent formula at
\(R-2\) and then the second invariant for \(P_{t-1}\) gives
\[
 Z_{b-1,b}=E_bZ_{b-1,b+1}-v^{-1}Z_{b-1,b+1}E_b.
\]
Starting from \(Z_{R-1,R}=E_R\), these two recurrences give
\eqref{eq:endpoint-root-recursion}.  The endpoint-raising operators
\(D_a\) commute with \(D_b,\ldots,D_{R-1}\) whenever
\(a\leq b-2\), so the blockwise recurrences have exactly the
operator order displayed there.

There is also a representation-theoretic way to check all endpoint
coefficients, including the lowering coefficients.  For the opposite
adjoint action, \(E_R\) is a lowest Levi-weight vector of weight
\(\varepsilon_{R-1}+\varepsilon_R\).  Every \(F_i\), \(i\in J_R\),
kills it.  Every \(E_i\) except \(E_{R-2}\) kills it, and
\(E_{R-2}^2\) kills it by the quantum Serre relation
\eqref{eq:quantum-serre-relations}.  These are the defining
lowest-weight relations of the minuscule Levi module
\(\bigwedge_v^2 V_R\).  The submodule generated by \(E_R\) is
therefore that simple module.  Recursion
\eqref{eq:endpoint-root-recursion} gives all its weight vectors with
coefficient one.  The rank-one relation on each two-element string gives
the reverse \(F_i\)-coefficient one.  This proves the endpoint formulas.

The opposite adjoint action has the module-algebra identities
\[
 \begin{aligned}
 E_i(ab)&=(E_i a)b+(K_i a)(E_i b),\\
 F_i(ab)&=a(F_i b)+(F_i a)(K_i^{-1}b).
 \end{aligned}
\]
Since it preserves the span of the generators, it preserves their
algebra.  Reversing multiplication gives
\begin{equation}\label{eq:positive-star-leibniz}
 \begin{aligned}
 E_i(a\star b)&=(E_i a)\star(K_i b)+a\star(E_i b),\\
 F_i(a\star b)&=(F_i a)\star b+(K_i^{-1}a)\star(F_i b).
 \end{aligned}
\end{equation}
These are the product rules for \eqref{eq:coproduct-convention}.
On a Levi-trivial cyclic vector, the rightmost \(E_i,F_i\) terms in
\eqref{eq:opposite-adjoint-positive} vanish and \(K_i\) acts by one.
Consequently \((u\cdot a)v_0=u(av_0)\).  This proves agreement with
the parabolic quotient and with its sufficiently large spin embeddings.
\end{proof}

Let \(\mathscr A_d\) consist of the products supported on the first
\(d\) endpoints.  The endpoint formulas preserve this subspace under
the first \(d-1\) Levi generators.  PBW straightening and nonnegativity
of endpoint weights make it a subalgebra for \(\star\), and the based
quotient restricts to its homogeneous parts.  The spin approximation
and excellent restriction therefore apply with this explicit Levi action.

\begin{lemma}[Two straightening relations]
\label{lem:positive-two-relations}
In \(\mathscr A_d\),
\[
 Z_{bc}\star Z_{ac}=vZ_{ac}\star Z_{bc}\quad(a<b<c),\qquad
 Z_{ad}\star Z_{bc}=Z_{bc}\star Z_{ad}\quad(a<b<c<d).
\]
\end{lemma}
\begin{proof}
Apply the Levendorskii--Soibelman relation in the forward root order of
\(B_R^+\), then reverse the product.  For the first pair the root inner
product is one.  A correction would have weight
\(\varepsilon_a+\varepsilon_b+2\varepsilon_c\) and use only roots
strictly between \(\varepsilon_b+\varepsilon_c\) and
\(\varepsilon_a+\varepsilon_c\).  No such radical monomial exists:
both factors would have to contain \(c\), and their other endpoints
would have to be \(a,b\).  For the second pair the inner product is zero.
The two alternative pairings of its four endpoints each use a root
outside the intervening root interval.  Thus neither relation has a
correction term.  The leading LS factor gives the stated powers of \(v\).
\end{proof}

\subsubsection{The squarefree Hecke action}
\label{subsubsec:squarefree-hecke}

For $d=2p$, recall that $\delta_d=\sum_{a=1}^d\varepsilon_a$
(Subsection~\ref{subsec:squarefree-tools}).  At this weight the $p$
pairs use each endpoint exactly once.  The standard vector of the
resulting matching is
\begin{equation}\label{eq:matching-pbw}
 Z_w=\mathop{\star}_{(a,b)\text{ a pair of }w}Z_{ab},
 \qquad w\in\I_d^{\FPF},
\end{equation}
where the written factors have increasing right endpoints.  These are
the actual PBW representatives from
Lemma~\ref{lem:radical-pbw-balanced}, without additional scalar factors.

On squarefree weight, the \(i,i+1\) rank-one strings have highest
weight zero or two: the two endpoint coordinates have total two and
all endpoint weights are nonnegative.  With the normalization of
Subsection~\ref{subsec:squarefree-tools}, the Hecke operator is
\begin{equation}\label{eq:explicit-squarefree-hecke}
 \mathsf H_i=E_iF_i-v^{-1}\operatorname{id}
 \qquad\text{on squarefree weight}.
\end{equation}
It has eigenvalues \(-v^{-1}\) and \(v\) on those two zero-weight
lines.  In particular it preserves every Levi submodule.

\begin{lemma}[The complete matching calculation]
\label{lem:positive-matching-table}
The operators \eqref{eq:explicit-squarefree-hecke} have the
fixed-point-free column standard action \eqref{eq:four-cases}.
\end{lemma}
\begin{proof}
Put \(j=i+1\).  If \((i,j)\) is a pair, both \(E_i\) and \(F_i\)
kill that carrier, so \(\mathsf H_iZ_w=-v^{-1}Z_w\).
Otherwise there are two distinct carriers.  With \(A\) the lower
matching and \(B\) the matching obtained by swapping \(i,j\), their
affected factors and the weight-two target are
\[
\begin{array}{c|c|c|c}
 &A&B&Q\\\hline
 a<b<i&Z_{ai}\star Z_{bj}&Z_{bi}\star Z_{aj}&Z_{ai}\star Z_{bi}\\
 j<a<b&Z_{ia}\star Z_{jb}&Z_{ja}\star Z_{ib}&Z_{ia}\star Z_{ib}\\
 a<i<j<b&Z_{ai}\star Z_{jb}&Z_{aj}\star Z_{ib}&Z_{ai}\star Z_{ib}.
\end{array}
\]
In each row $A(i)=a<b=A(j)$, so conjugation by $s_i$
is a strict ascent in the descending completion.
The endpoint formulas, \eqref{eq:positive-star-leibniz}, and
Lemma~\ref{lem:positive-two-relations} give in all three cases
\[
 E_iA=Q,\qquad E_iB=vQ,\qquad F_iQ=B+v^{-1}A.
\]
For example, in the first row \(E_iB=Z_{bi}\star Z_{ai}=vQ\);
the first term of \(F_iQ\) is
\(Z_{aj}\star Z_{bi}=Z_{bi}\star Z_{aj}=B\).
In the other two rows the factor \(v\) in \(E_iB\) is the
\(K_i\)-weight of its second affected factor, and the two terms
of \(F_iQ\) are already in the displayed order.

These calculations apply inside a full matching.  Every other factor
is killed by \(E_i,F_i\) and has \(K_i\)-weight one.  It can therefore
be retained in its original position throughout each calculation.
When two affected right endpoints change their relative order, they are
\(i,j\), so no third right endpoint lies between them.
Thus the indicated straightening does not cross an unexamined factor.

The two zero-weight lines are spanned by
\(A-v^{-1}B\) and \(A+vB\), respectively.  Substitution in
\eqref{eq:explicit-squarefree-hecke} gives
\[
 \mathsf H_iA=B,\qquad
 \mathsf H_iB=A+(v-v^{-1})B.
\]
These and the pair case are the required standard table.
\end{proof}

\begin{lemma}[The minimum matching is bar fixed]
\label{lem:positive-minimum-bar}
The vector of \((1,2)(3,4)\cdots(d-1,d)\) in
\eqref{eq:matching-pbw} is bar invariant in \(Q_R\).
\end{lemma}
\begin{proof}
By Lemma~\ref{lem:radical-pbw-balanced}, each surviving correction
has the same weight and strictly larger radical data in left
lexicographic order.  At squarefree weight these data are matchings.

The indicated matching is lexicographically maximal.  Its first
available radical root is \(\varepsilon_{d-1}+\varepsilon_d\),
whose exponent is one.  Any other matching omitting this pair has
exponent zero at the first possible difference.  Among matchings
containing it, delete the pair and apply the same argument to the
first \(d-2\) endpoints.  This induction proves maximality.
Coordinates involving endpoints beyond \(d\) are zero throughout.
There is therefore no surviving nontrivial bar correction.
\end{proof}

\subsubsection{Singletons and the based tensor product}

We use the lower global basis of \(V(k\omega_1)\), with highest vector
\(v_k\).  For a composition \(\mathbf a=(a_1,\ldots,a_d)\) of \(k\),
write
\[
 s_{\mathbf a}=
 F_{d-1}^{(a_d)}F_{d-2}^{(a_{d-1}+a_d)}\cdots
 F_1^{(a_2+\cdots+a_d)}v_k.
\]
Its weight spaces are one dimensional; these are the lower global
vectors obtained by the indicated divided-power strings.  Their
\(\A\)-span and \(\Oinf\)-span are the integral and lower crystal
lattices.  Their normalization gives
\[
 E_i s_{\mathbf a}=[a_i+1]s_{\mathbf a+\epsilon_i-\epsilon_{i+1}},
 \qquad
 F_i s_{\mathbf a}=[a_{i+1}+1]
 s_{\mathbf a-\epsilon_i+\epsilon_{i+1}},
\]
when the source coordinate being removed is positive, and zero
otherwise.  For a set \(F\) of size \(k\), let \(s_F\) denote
the vector whose composition is its indicator function.

Put \(d=2p+k\) and
\begin{equation}\label{eq:tensor-model}
 \mathscr X_{d;p,k}=\mathscr A_d^{(p)}\otimes V(k\omega_1).
\end{equation}
We use the coproduct \eqref{eq:coproduct-convention} and the lower
based tensor-product bar \eqref{eq:tensor-product-bar}
\cite{Lusztig1992Tensor,Lusztig1993}.
As in Subsection~\ref{subsec:squarefree-tools}, squarefreeness is a
condition on the total weight.  Each label occurs exactly once across
the pair and singleton factors.  The standard basis
at weight $\delta_d$ is
\begin{equation}\label{eq:partial-matching-tensor}
 Z_w=\left(\mathop{\star}_{(a,b)\text{ a pair of }w}Z_{ab}\right)
 \otimes s_{\Fix(w)},
 \qquad w\in\I_d,\quad\#\Fix(w)=k.
\end{equation}
The first tensor factor uses precisely the moving endpoints and the
second precisely the fixed endpoints.  The tensor-product crystal
lattice is the tensor product of the two lower lattices.  Hence these
vectors are an integral basis and a lower-lattice basis of the slice,
and their residues are its crystal basis.

\begin{theorem}[The type-$D$ column based comparison]
\label{lem:based-block}
The weight space $(\mathscr X_{d;p,k})_{\delta_d}$, with
\(\mathsf H_i=E_iF_i-v^{-1}\), its integral lattice, and its lower
global basis, is based-isomorphic to the fixed-\(k\) column block.
\end{theorem}
\begin{proof}
For pair--pair carriers use
Lemma~\ref{lem:positive-matching-table}.  For pair--singleton carriers
write
\[
 A=x_i\otimes s_{F\cup\{i+1\}},\qquad
 B=x_{i+1}\otimes s_{F\cup\{i\}},\qquad
 Q=x_i\otimes s_{F\cup\{i\}},
\]
where the other endpoint of the affected pair and all unaffected
factors are retained, and \(F\) contains neither \(i\) nor \(i+1\).
In the descending completion, the partner of $i$ in $A$ lies
among the original labels, whereas the partner of the singleton $j=i+1$
is an added label.  Thus this is the strict-ascent case of
\eqref{eq:four-cases}.
The coproduct and the displayed symmetric-power coefficients give
\(E_iA=Q\), \(E_iB=vQ\), and
\(F_iQ=B+v^{-1}A\).  Thus the same two-by-two Hecke matrix results.
If both endpoints are singletons, the local zero-weight line is the
middle of the symmetric \(V(2)\) string, so its eigenvalue is \(v\).
Together with the pair eigenvalue \(-v^{-1}\), these are exactly
the four cases in \eqref{eq:four-cases}.  In particular, the Hecke
relations hold because this is the prescribed column standard action.

Write \(r=2p\).  The minimum standard tensor is
\[
 (Z_{12}\star Z_{34}\star\cdots\star Z_{r-1,r})
 \otimes s_{\{r+1,\ldots,d\}}.
\]
Its factors are bar fixed, the first by
Lemma~\ref{lem:positive-minimum-bar} and the second by its definition
as a global vector.  By \eqref{eq:tensor-product-bar}, a term
$\Theta_\beta$ shifts the two weights by $\beta$ and $-\beta$.
Write such a weight as $\beta=\sum c_i\alpha_i$, with
$c_i\geq0$ and $\beta\ne0$, and let $s$ be the largest index with
$c_s>0$.  If $s\geq r$, applying the positive-weight term to the
first factor would give endpoint coordinate $-c_s$ at $s+1$.
This is impossible in its polynomial weight set.  If $s<r$, let
$t$ be the smallest index with $c_t>0$.  Applying the negative-weight
term to the singleton factor would give coordinate $-c_t$ at $t$,
which is again impossible.  Thus every nonconstant term vanishes,
and the tensor bar fixes the minimum.

The module bar fixes \(E_i,F_i\), so
\[
 \overline{\mathsf H_i x}=(E_iF_i-v)\bar x
 =\mathsf H_i^{-1}\bar x.
\]
Both this bar and the column-module bar are compatible with the
same standard action and fix its minimum.  The strict-ascent
argument in Subsection~\ref{subsec:reduced-graphs} therefore identifies
the two bars.  The standard integral lattices and lower crystal
lattices, including their residues, were identified in
\eqref{eq:partial-matching-tensor}.
By \eqref{eq:canonical-normalization}, the column canonical vectors
are bar-fixed lifts of these residues.  Balanced uniqueness
\eqref{eq:balanced-uniqueness-prelim} identifies them with the lower
global vectors, proving the based identification.
\end{proof}

\subsection{Verification of the three conditions}

\begin{lemma}[Multiplicity-free tensor character]
\label{lem:tensor-character}
The type-$A_{d-1}$ character of the tensor model is
\begin{equation}\label{eq:pieri-shapes}
 \operatorname{ch}\mathscr X_{d;p,k}
 =\sum_{\substack{\mu\vdash d\\\mu\text{ has exactly }k\text{ odd columns}}}s_\mu(X),
\end{equation}
and every summand occurs with multiplicity one.
\end{lemma}
\begin{proof}
Give each radical generator degree one, recorded by $t$.  The root
vector $Z_{ij}$ has weight $\varepsilon_i+\varepsilon_j$.  Summing over
its exponent in a divided PBW monomial gives the factor
$(1-tX_iX_j)^{-1}$.  The exponents vary independently, so the graded
character is the product on the left of Littlewood's identity
\[
 \prod_{i<j}(1-tX_iX_j)^{-1}
 =\sum_{\nu:\,\nu^{\tr}\text{ has even parts}}
   t^{|\nu|/2}s_\nu(X).
\]
The expansion follows by induction on the number of variables
\cite[Chapter~I]{Macdonald1995}.  After adjoining $X_d$, expand the new factors as complete
symmetric functions.  The coefficient of
$s_\rho(X_{<d})X_d^q$ is nonzero precisely when $\rho/\nu$ is a horizontal
$q$-strip and $\nu$ has even columns.  Such $\nu$ is unique: delete the
bottom box of every odd column of $\rho$.  If $q$ is the number of odd
columns of $\rho$, the coefficient is $t^{(|\rho|+q)/2}$.  The branching
identity
\[
 s_\lambda(X_1,\ldots,X_d)
 =\sum_{\rho:\,\lambda/\rho\text{ is horizontal}}
   s_\rho(X_{<d})X_d^{|\lambda|-|\rho|}
\]
has the same coefficients, since the unique even-column partition above
$\rho$ is obtained by adding a bottom box to every odd column.  This proves
the product identity and its multiplicity-one assertion.

Multiplication by
$\operatorname{ch}\operatorname{Sym}^k(V_d)=h_k=s_{(k)}$, followed by the
horizontal-strip rule and extraction of the coefficient of $t^p$, gives
\eqref{eq:pieri-shapes}.  For a partition $\mu$ on the right, deleting the
bottom box of every odd column gives its unique even-column predecessor;
conversely, the deleted boxes form a horizontal strip.  Thus every summand
occurs exactly once.
\end{proof}

For a column vertex $w$ in the fixed-$k$ block, let $\topJ(w)$ be
its lower Levi top in $\mathscr X_{d;p,k}$ under
Theorem~\ref{lem:based-block}.  Put
\[
 \Lambda^{\mathrm c}_{d;k}
 =\{\mu\vdash d:\mu\text{ has exactly }k\text{ odd columns}\}.
\]

\begin{proposition}[The three conditions for the column model]
\label{prop:column-criterion-hypotheses}
Fix $d=2p+k\geq2$.  The fixed-$k$ column block satisfies conditions
\ref{cond:levi-top-multiplicity}--\ref{cond:levi-top-connectivity}
for $\tau(w)=\topJ(w)$ and $\Lambda=\Lambda^{\mathrm c}_{d;k}$,
ordered by dominance.  For every dominance upset
$\Omega\subseteq\Lambda^{\mathrm c}_{d;k}$, the space
\begin{equation}\label{eq:column-positive-top-upset}
 \mathscr M^{\mathrm c}_\Omega
 =\bigoplus_{\topJ(w)\in\Omega}\A\,\underline N_w
\end{equation}
is a saturated $\Hh_d$-submodule.  For each occurring top $\mu$,
the one-top subquotient is based-isomorphic to the squarefree lower
slice $V(\mu)_{\delta_d}$, and its top fiber is bidirected connected
in the original column graph.
\end{proposition}
\begin{proof}
For condition~\ref{cond:levi-top-multiplicity},
Lemma~\ref{lem:tensor-character} gives precisely the highest weights
in $\Lambda^{\mathrm c}_{d;k}$, each with multiplicity one.
For each such $\mu$, recall the
\hyperref[not:squarefree-tableaux]{tableau indexing of the lower global
basis}.  Its weight-$\delta_d$ vectors are the $G_S$ with
$S\in\SYT(\mu)$ by Lemma~\ref{lem:tableau-slice}, so the slice is
nonempty.
Together with Theorem~\ref{lem:based-block}, this establishes the
surjective multiplicity-free top datum.

For condition~\ref{cond:levi-top-upsets},
Theorem~\ref{thm:excellent} gives the saturated Levi submodule of
$\mathscr X_{d;p,k}$ spanned by tops in $\Omega$.
The action $\mathsf H_i=E_iF_i-v^{-1}\operatorname{id}$ in
Theorem~\ref{lem:based-block} preserves that submodule and its
squarefree slice.  The based comparison identifies the resulting
slice with exactly \eqref{eq:column-positive-top-upset}.
Its quotient is free on the complementary column canonical vectors,
so it is saturated.

For condition~\ref{cond:levi-top-connectivity}, fix
$\mu\in\Lambda^{\mathrm c}_{d;k}$ and set
\[
 \Omega_1=\{\nu\in\Lambda^{\mathrm c}_{d;k}:\mu\preceq\nu\},
 \qquad
 \Omega_0=\Omega_1\setminus\{\mu\}.
\]
Excellent restriction and multiplicity one identify
$\mathscr M^{\mathrm c}_{\Omega_1}/
 \mathscr M^{\mathrm c}_{\Omega_0}$ with the squarefree lower based
slice of the unique $V(\mu)$.
Corollary~\ref{cor:lower-layer-connected} gives bidirected
connectivity in that slice.  Since the filtration uses actual
canonical-basis subsets, the coefficients between surviving
vectors are unchanged in the quotient.  The bidirected edges
therefore lift to the original column graph, proving connectivity
of its whole top fiber.
\end{proof}

\subsection{Cells, insertion shapes, and dominance}
\label{subsec:column-results}

Proposition~\ref{prop:column-criterion-hypotheses} verifies
conditions~\ref{cond:levi-top-multiplicity}--\ref{cond:levi-top-connectivity}
for the column realization constructed in this section.  We now
apply Theorem~\ref{thm:levi-top-criterion} and identify its top fibers
by column insertion.  As in the row model, the conclusions include
$d=1$, when the graph has one vertex with polynomial top $(1)$
and every assertion is immediate.

\begin{theorem}[$\ncell$-cells, $\ncell$-molecules, and column Levi tops]
\label{thm:n-cells-molecules}
In each block with $k$ fixed points, the $\ncell$-cells and
$\ncell$-molecules are precisely the lower Levi-top fibers
\[
 C^{\mathrm c}_\mu
 =\{w\in I_d:\#\Fix(w)=k,\ \topJ(w)=\mu\},
\]
where $\mu$ ranges over the occurring tops.  In particular, every
$\ncell$-cell of $\nGamma$ is exactly one $\ncell$-molecule.
\end{theorem}
\begin{proof}
For $d\geq2$, Proposition~\ref{prop:column-criterion-hypotheses}
verifies all three hypotheses of Theorem~\ref{thm:levi-top-criterion}
in each fixed-$k$ block.  The criterion therefore identifies the
cells, molecules, and the
fibers $C^{\mathrm c}_\mu$.  By
Subsection~\ref{subsec:reduced-graphs}, no coefficient arrow crosses
fixed-point blocks, so the assertion also holds in the full column graph.
\end{proof}

\subsubsection{The column Levi label}

Fix $k$, put $p=(d-k)/2$, and work in the fixed-$k$ block.  For a top
$\mu$ in \eqref{eq:pieri-shapes}, let $C^{\mathrm c}_\mu$ be the fiber
defined in Theorem~\ref{thm:n-cells-molecules}.
By Proposition~\ref{prop:column-criterion-hypotheses}, its one-top
associated graded is the squarefree lower based slice of $V(\mu)$.
Recall that $\nAsc(w)$ records the indices $i$ with
\hyperref[not:squarefree-tableaux]{diagonal coefficient $-v^{-1}$}.
This coefficient is unchanged on passing from $\underline N_w$ to its
image $G_S$ in the layer.  By Lemma~\ref{lem:tableau-slice} those
indices are precisely the descents of $S$.  Hence
\eqref{eq:schur-fundamental} gives
\begin{equation}\label{eq:top-character}
 \sum_{w\in C^\mathrm c_\mu}F_{\nAsc(w)}
 =\sum_{S\in\SYT(\mu)}F_{\Des(S)}=s_\mu.
\end{equation}

By Theorem~\ref{thm:n-cells-molecules}, each
$C^{\mathrm c}_\mu$ is exactly one column molecule.
Proposition~\ref{prop:reduced-molecules} identifies that molecule
with a single column insertion-shape fiber, say of shape $\lambda$.
Comparing \eqref{eq:top-character} with the column identity
\eqref{eq:column-shape-character-new} in
Proposition~\ref{prop:insertion-shape-characteristics} gives
$s_\mu=s_{\lambda^{\tr}}$.
The Schur functions are linearly independent, so
\begin{equation}\label{eq:label-identification}
 \topJ(w)=\lambda_{\mathrm c}(w)^{\tr}.
\end{equation}

This also verifies that the excellent filtration contains exactly one copy of
each layer predicted by \eqref{eq:pieri-shapes}.

\subsubsection{Column dominance and the cells}

We can now state the filtration in the original tableau notation.

\begin{theorem}[Column dominance filtration]\label{thm:dominance-filtration}
Let $\mathscr D$
be a downset in the dominance order on partitions of $d$.
Then
\begin{equation}\label{eq:downset-module}
 \mathcal F_{\mathscr D}
 =\bigoplus_{\lambda_{\mathrm c}(w)\in\mathscr D}\A\,\underline N_w
\end{equation}
is a saturated $\Hh_d$-submodule of the column model.
\end{theorem}
\begin{proof}
For each admissible $k$, put
\[
 \mathcal F^{\mathrm c}_{\mathscr D,k}=
 \bigoplus_{\substack{\#\Fix(w)=k\\
             \lambda_{\mathrm c}(w)\in\mathscr D}}
             \A\,\underline N_w.
\]
By \eqref{eq:label-identification}, its Levi-top set is
$\Omega_k=\{\mu:\mu^{\tr}\in\mathscr D\}$ among the tops occurring in
the fixed-$k$ block.  Equation~\eqref{eq:transpose-reverses} makes
$\Omega_k$ an upset.  Proposition~\ref{prop:column-criterion-hypotheses} makes the
displayed canonical-basis span stable and saturated.  Therefore
\[
 \mathcal F^{\mathrm c}_{\mathscr D}
 =\bigoplus_{k\equiv d\pmod2}\mathcal F^{\mathrm c}_{\mathscr D,k}
\]
is a stable saturated submodule of the complete column model; this is
exactly \eqref{eq:downset-module}.
\end{proof}

\begin{theorem}[Column edge dominance]\label{thm:column-edge-dominance}
Let $y,z\in\I_d$ and put
$\widetilde y=\iota_{\mathrm{des}}(y)$ and
$\widetilde z=\iota_{\mathrm{des}}(z)$.  If
$\widetilde y\longrightarrow_{\ncell}\widetilde z$ is a coefficient arrow
in the reduced column Gelfand $W$-graph, then
\[
 \lambda_{\mathrm c}(z)\preceq\lambda_{\mathrm c}(y).
\]
If the reverse arrow is also present, applying the same implication in both
directions gives equality of the two shapes.
\end{theorem}
\begin{proof}
Let $\alpha=\lambda_{\mathrm c}(y)$ and take the principal downset
\[
 \mathscr D_\alpha=\{\beta:\beta\preceq\alpha\}.
\]
The vector $\underline N_y$ belongs to the stable submodule
$\mathcal F_{\mathscr D_\alpha}$ of
Theorem~\ref{thm:dominance-filtration}.  By the definition of the arrow,
$\underline N_z$ has a nonzero coefficient in
$H_i\underline N_y$ for some $i$.  Since
$\mathcal F_{\mathscr D_\alpha}$ is spanned by a subset of the canonical
basis, that coefficient can be nonzero only if
$\lambda_{\mathrm c}(z)\in\mathscr D_\alpha$.  Hence
$\lambda_{\mathrm c}(z)\preceq\lambda_{\mathrm c}(y)$.
\end{proof}

\begin{corollary}[Column cells by insertion shape]\label{cor:column-cells-shape}
Two vertices lie
in the same $\ncell$-cell if and only if their modified column
Beissinger tableaux have the same shape.  Equivalently, every $\ncell$-cell is
one $\ncell$-molecule.
\end{corollary}
\begin{proof}
Theorem~\ref{thm:n-cells-molecules} identifies cells with molecules.
Proposition~\ref{prop:reduced-molecules} identifies the molecules
of this reduced graph with the fibers of the modified column insertion shape.
\end{proof}

\begin{proof}[Completion of the proof of Theorem~\ref{thm:main-intro}]
Corollary~\ref{cor:column-cells-shape} and
Theorem~\ref{thm:column-edge-dominance} prove part~\textup{(b)}.
Part~\textup{(a)} was proved in Subsection~\ref{subsec:row-results}.
Theorems~\ref{thm:row-dominance-filtration}
and~\ref{thm:dominance-filtration} give the two saturated
canonical-basis filtrations asserted in the final sentence.
\end{proof}

\newpage
\subsection{A rank-six column cell}

\begin{example}[A complete $\ncell$-cell $W$-graph in rank six]
\label{ex:rank-six-n-cell}
Work in the block of $I_6$ with two fixed points.  Define
\[
\begin{array}{c|c|c|c}
 x&w_x&\PcB(w_x)&\nAsc(w_x)\\ \hline
 A&(13)(24)(5)(6)&\rowtableau{1&2&5\cr3&4&6\cr}&\{1,3,4\}\\
 B&(13)(25)(4)(6)&\rowtableau{1&2&4\cr3&5&6\cr}&\{1,3,5\}\\
 C&(14)(23)(5)(6)&\rowtableau{1&3&5\cr2&4&6\cr}&\{2,4\}\\
 D&(14)(25)(3)(6)&\rowtableau{1&2&3\cr4&5&6\cr}&\{1,2,4,5\}\\
 E&(15)(23)(4)(6)&\rowtableau{1&3&4\cr2&5&6\cr}&\{2,3,5\}.
\end{array}
\]
These are exactly the five vertices with restricted column insertion shape
$(3,3)$.  Hence
Proposition~\ref{prop:reduced-molecules} and
Theorem~\ref{thm:n-cells-molecules} identify
\[
 \mathcal C_{\ncell}^{(3,3)}
 =\{w_A,w_B,w_C,w_D,w_E\}
\]
as one $\ncell$-cell and one $\ncell$-molecule.

Applying the standard action \eqref{eq:four-cases} to the column
canonical basis normalized by \eqref{eq:canonical-normalization}
gives the graph below.  Solid double arrows
are bidirected, and the dashed arrow is one-directional.  The directed-edge table lists all
generators producing each coefficient.  Every displayed coefficient has
weight $1$, and there are no further off-diagonal arrows.  Thus the two
tables and the diagram give the complete reduced cell-subquotient $W$-graph.
\[
\begin{tikzpicture}[
  x=1.25cm,y=1.05cm,
  vertex/.style={circle,draw,thick,minimum size=7mm,inner sep=0pt},
  bi/.style={{Stealth[length=2.1mm]}-{Stealth[length=2.1mm]},semithick},
  one/.style={-{Stealth[length=2.1mm]},thick,dashed}
]
 \node[vertex] (A) at (0,2.5) {$A$};
 \node[vertex] (B) at (-2,1) {$B$};
 \node[vertex] (C) at (2,1) {$C$};
 \node[vertex] (E) at (0,-0.5) {$E$};
 \node[vertex] (D) at (-2,-2.2) {$D$};

 \draw[bi] (A) -- (B);
 \draw[bi] (A) -- (C);
 \draw[bi] (B) -- (E);
 \draw[bi] (C) -- (E);
 \draw[bi] (B) -- (D);
 \draw[one,bend left=35] (D) to (C);
\end{tikzpicture}
\]
\[
\begin{array}{c|c@{\qquad}c|c@{\qquad}c|c}
\text{arrow}&J&\text{arrow}&J&\text{arrow}&J\\ \hline
A\to B&\{4\}&B\to A&\{5\}&A\to C&\{1,3\}\\
C\to A&\{2\}&B\to E&\{1\}&E\to B&\{2\}\\
C\to E&\{4\}&E\to C&\{3,5\}&B\to D&\{3\}\\
D\to B&\{2,4\}&D\to C&\{1,5\}&\multicolumn{2}{c}{}
\end{array}
\]
The bidirected part is connected, and the curved arrow $D\to C$ is the
unique additional one-directional edge in this five-vertex cell graph.
\end{example}

\end{document}